\documentclass{amsart}
\usepackage{amssymb}
\usepackage{amsfonts}
\usepackage{amsmath}
\usepackage{lastpage}
\usepackage[legalpaper,bookmarks=true,colorlinks=true,linkcolor=blue,citecolor=blue]{hyperref}
\usepackage{graphicx}%
\usepackage{fancyhdr}
\usepackage{color}
\usepackage[mathlines]{lineno}
\usepackage{lscape}
\usepackage{epsfig}
\usepackage{natbib}
\usepackage{geometry}
\usepackage{tgbonum}
\fontfamily{qcr}\selectfont
\usepackage[]{listings}
\usepackage[normalem]{ ulem }
\usepackage{soul}

\newcommand{\NC}{ {\ \neg \mathcal{R} \ }}

\newcommand{\Bin}{\bigskip \noindent}
\newcommand{\Bi}{\bigskip}
\newcommand{\Ni}{\noindent}

\newtheorem{theorem}{Theorem}
\theoremstyle{plain}

\newtheorem{corollary}{Corollary}

\newtheorem{lemma}{Lemma}

\newtheorem{proposition}{Proposition}

\numberwithin{equation}{section}

\begin{document}
\Large
\title[Introduction to General Theory of records]{An introduction to a general records theory both for dependent data and high dimensions, revisited}
\author{Gane Samb Lo}
\author{Mohammad Ahsanullah}

\begin{abstract} .\\  
The probabilistic investigation on record values and record times of a sequence of random variables defined on the same probability space has received much attention from 1952 to now. A great deal of such theory focused on \textit{iid} or independent real-valued random variables. There exists a few results for real-valued dependent random variables. Some papers deal also with multivariate random variables. But a large theory regarding vectors and dependent data has yet to be done. In preparation of that, the probability laws of records are investigated here, without any assumption on the dependence structure. The results are extended sequences with values in partially ordered spaces whose order is compatible with measurability. The general characterizations are checked in known cases mostly for \textit{iid} sequences. The frame is ready for undertaking a vast study of records theory in high dimensions and for types of dependence.\\

\noindent $^{\dag}$ Gane Samb Lo.\\
LERSTAD, Gaston Berger University, Saint-Louis, S\'en\'egal (main affiliation).\newline
LSTA, Pierre and Marie Curie University, Paris VI, France.\newline
AUST - African University of Sciences and Technology, Abuja, Nigeria\\
gane-samb.lo@edu.ugb.sn, gslo@aust.edu.ng, ganesamblo@ganesamblo.net\\
Permanent address : 1178 Evanston Dr NW T3P 0J9,Calgary, Alberta, Canada.\\

\noindent $^{\dag\dag}$ Mohammad Ahsanullah\\
Department of Management Sciences. Rider University. Lawrenceville, New Jersey, USA\\
Email : ahsan@rider.edu\\


\textcolor{white}{
a\\
a\\
a\\
a\\
a\\
a\\
a\\
a\\
a\\
a\\
a\\
a\\
a\\
a\\
a\\
}
\end{abstract}
\maketitle

\newpage
\tableofcontents

\newpage
\noindent \textbf{R\'esum\'e}. L'investigation probabiliste des records d'une suite de variables aleatoires a reçu une grande attention depuis 1952 jusqu'à nos jours. Une grande partie de cette théorie a concerné les suites de variables aléatoires indépendantes, identiquement distribuées ou non, à valeurs réelles. Il existe quelques résultats pour les variables aléatoires dépendantes à valeurs réelles. Certains articles traitent également de variables aléatoires multivariées. Cependant, une théorie majeure relative aux données dépendantes ou multivariées est encore à faire. En préparation de cela, les lois de probabilité des records sont étudiées ici, sans aucune hypothèse sur la structure de dépendance des variables. Les résultats sont étendus aux suite de variables à valeurs dans un espace partiellement ordonné sur lequel l'ordre est compatible avec la mesurabilité. Les caractérisations générales sont testées sur les résultats connus portant principalementent sur les suites \textit{iid}. Le cadre est prêt pour passer à la théorie en haute dimension et pour differents types de dépendance.\\

\section{Introduction}

\noindent The theory of records both dealing with record values and record times for a sequence of random variables $(X_n)_{n\geq 1}$, defined on the same probability space $(\Omega, \mathcal{A}, \mathbb{P})$ and taking their values in some measurable space $E$ endowed with a partial order ($\leq$) is a relatively recent sub-discipline of Probability Theory and Statistics. That theory goes back to \cite{chandler} and back to \cite{feller} who applied it to gambling problems.\\

\noindent By now, the theory has had extraordinary developments in a great variety of directions including characterization of distributions or of stochastic processes, statistical estimations. A few number of attempts beyond the scheme of independent and identically distributed (\textit{iid}) sequences are available.\\

\noindent A stochastic process view has led to the extremal process (see \cite{dwass}) which has helped to solve many problems in Extreme Value Theory (see \cite{resnick87}).\\

\noindent Especially, in the \textit{iid} case, that theory is tremendously developed in a number of papers and more than a dozen of books (as quoted by \cite{ahsanullahUGB}) has been reported. To cite a few, we have the following books : \cite{ahsan08}, \cite{ahsan88}, \cite{ahsan95}, \cite{ahsan04}, \cite{ahsan06}, \cite{arnold}, \cite{gulatis}, \cite{nevzorov}, \cite{resnick87} (partially).\\

\noindent Although several departures from the real-valued sequences frame, from the \textit{iid} case and even from the independent assumption, it seems that there does not exist a complete review on the results of such generalizations, including the ones related to partial order relations, for examples in finite dimensional spaces like $\mathbb{R}^d$, $d>1$.\\

\noindent The aim of this paper is two-fold. First, we wish to propose a general frame for the theory of records values and record times in an arbitrary partially or totally ordered space, including random fields and to find out the finite-dimensional probability laws for the records values and times. Such a general frame would be an appropriate place to summarize available generalizations and to be the basis to get new extensions.\\


\noindent It is expected that this presentations and the general formulas therein  will quick off new trends of innovative research works from the readers.\\

\section{General setting} \label{gsetting}


\subsection{Basic definitions about records on $\mathbb{R}$} \label{gs_def}

\noindent We are going to introduce all the needed definitions on a sequence of real numbers $x=(x_{n})_{n\geq 1}$.\\

\noindent \textbf{(A) - Strong record times and strong upper records}.\\

\noindent Let us  define them by induction. The first record time, in general, is set to one, and we write $u(1)=1$ and the first record is defined by

\begin{equation*}
x^{(1)}=x_{u(1)}=:x_{1}.
\end{equation*}

\Bin Next, we search

\begin{equation*}
u(2)=\left\{ 
\begin{tabular}{ll}
$\inf \{j>u(1),x_{j}>x_{u(1)}\}\equiv \inf A_2$ &  if $A_2\neq \emptyset$\\
$+\infty $& otherwise \\ 
\end{tabular}
\right. .
\end{equation*}

\noindent If $u(2)<+\infty$, we call $u(2)$ the second strong record time and the strong record value is given by

\begin{equation*}
x^{(2)}=x_{u(2)}.
\end{equation*}

\bigskip
\noindent Given that the $n$-th strong record time exists, we may define

\begin{equation*}
u(n+1)=\left\{ 
\begin{tabular}{lll}
$\inf \{j>u(n),x_{j}>x_{u(n)}\}\equiv A_n$ & if $A_n \neq \emptyset$ \\
=$+\infty $ & otherwise \\ 
\end{tabular}
\right. .
\end{equation*}

\bigskip \noindent And, as previously, $u(n+1)$ is the $(n+1)$th strong record time if $u(n+1)<+\infty$, and 

\begin{equation*}
x^{(n+1)}=x_{u(n+1)}
\end{equation*}

\bigskip \noindent is the $(n+1)-th$ strong record value.\\

\noindent Either we may proceed indefinitely and the sequence of record times $(u(n))_{n\geq 1}$ is unbounded or we stop the
first time we have $u(n)=+\infty $ and in such a case, the sequence $x$ has only $(n-1)$ record  values where $n$ is necessarily greater than  one.\\

\bigskip \noindent \textbf{(B) - Weak records times and weak upper record values}.\\

\noindent We may define weak versions of record times and values by allowing repetitions of the record values. The first weak record time set to one, and we write $u^{(w)}(1)=1$ and the weak first record value is defined by

\begin{equation*}
x^{(1,w)}=x_{u^{(w)}(1)}=:x_{1}.
\end{equation*}

\Bin Next we define

\begin{equation*}
u(2)=\left\{ 
\begin{tabular}{ll}
$\inf  \{j>u^{(w)}(1),x_{j}\geq x_{u(1)}\}\equiv \inf A_{2} $ & if $A_{n-1} \neq \emptyset$  \\
$+\infty $ & otherwise\\
\end{tabular}
\right. .
\end{equation*}

\Bin If $u^{(w)}(2)<+\infty $, we call $u^{(w)}(2)$ the second weak record

time and the weak record  value is given by
\begin{equation*}
x^{(2,w)}=x_{u^{(w)}(2)}.
\end{equation*}

\bigskip
\noindent Given that the $n$-th weak record time exists, we may define the 
\begin{equation*}
u^{(w)}(n+1)=\left\{ 
\begin{tabular}{ll}
$\inf \{j>u^{(w)}(n),x_{j}\geq x_{u(n)}\}=\inf A_n$  &  if $A_n\neq \emptyset$\\
$+\infty $ & otherwise\\ 
\end{tabular}
\right. .
\end{equation*}

\bigskip  \noindent Next we have, as previously, that if $u^{(w)}(n+1)$ is the $(n+1)-th$ weak record time if $u^{(w)}(n+1)<+\infty $, we may define 

\begin{equation*}
x^{(n+1,w)}=x_{u^{(w)}(n+1)}
\end{equation*}

\noindent is the $(n+1)$-th weak strong record  value.\\

\bigskip \noindent \textbf{C - Lower records}.\\

\noindent Generally we set the first record time to one for any kind of record  value. The strong lower record times $(\ell (n))_{n\geq 1}$ and lower record values $(x_{\ell (n)})_{n\geq 1}$ are similarly defined by setting $\ell (1)=1$ and $x_{(1)}=x_{\ell (1)}=1$ and next by induction, if $\ell (n)$ is finite, by 

\begin{equation*}
\ell (n+1)=\left\{ 
\begin{tabular}{ll}
$\inf \{j>\ell (n),x_{j}<x_{\ell (n)}\}\equiv A_n$ &  if $A_2 \neq \emptyset$\\
$+\infty $ & otherwise 
\end{tabular}
\right. .,
\end{equation*}

\noindent and $\ell (n+1)$ is the $(n+1)$-th strong lower record time if it is finite
and the $(n+1)-th$ is given by 

\begin{equation*}
x_{(n+1)}=x_{\ell (n+1)}.
\end{equation*}

\noindent As well,  the weak lower record times $(\ell ^{(w)}(n))_{n\geq 1}$ and lower record  value $(x_{\ell ^{(w)}(n)})_{n\geq 1}$ are similarly defined by setting $\ell ^{(w)}(1)=1$ and $y_{(n,w)}=x_{\ell ^{(w)}(1)}=1$ and next by induction, if $\ell^{(w)}(n)$ is finite, by

\begin{equation*}
\ell ^{(w)}(n+1)=\left\{ 
\begin{tabular}{lll}
$+\infty $ & $if$ & $\{j>\ell ^{(w)}(n),x_{j}\leq x_{\ell^{(w)}(n)}\}=\emptyset $ \\ 
$\inf  \{j>\ell ^{(w)}(n),x_{j}\leq x_{\ell ^{(w)}(n)}\}$ &  & otherwise.%
\end{tabular}
\right. ,
\end{equation*}

\bigskip \noindent and $\ell ^{(w)}(n+1)$ is the $(n+1)$-th weak lower record time if it is
finite and the $(n+1)$-th lower record value is given by 
\begin{equation*}
x_{(n+1,w)}=x_{\ell ^{(w)}(n+1)}.
\end{equation*}

\bigskip
\subsection{Random records and general formulas}

\noindent We are going to move from statistic records to random ones.\\\

\noindent \textbf{(I) - Introduction}.\\

\noindent To make short, we call \textit{record values} by \textit{records} simply and \textit{records}, without any further precision, are \textit{strong upper records}.\\

\noindent In this section, we consider now a sequence of random real random variables

\begin{equation*}
X_{1}, \ X_{2},\ldots
\end{equation*}

\noindent defined on the same probability space $(\Omega ,\mathcal{A},\mathbb{P})$. We then define the random records associated to this sequence. The records
times and records become random variables defined by capital letters as follows
\begin{eqnarray*}
&&U(n), \ L(n),\text{ }U^{(w)}(n), \ L^{(w)}(n),\\
&&X^{(n)}=X_{U(n)},\text{ }Y^{(n)}=X_{L(n)},\\
&&X^{(n,w)}=X_{U^{(w)}(n)},\\
&&X_{(n,w)}=X_{L^{(w)}(n)}\\
\end{eqnarray*}

\noindent We are going to give a series of general facts on the laws of the record values and record times. We may focus only of the upper records, since we are able to derive the results on lower records from those on upper records by using the transform 
$$
(X_{1}, \ X_{2},\ldots) \ \ \mapsto \ \ (-X_{1}, \ -X_{2},\ldots)
$$

\noindent or, if the $X_i$'s are \textit{a.s.} positive, by using the transform 

$$
(X_{1}, \ X_{2},\ldots) \ \ \mapsto \ \ (1/X_{1}, \ 1/X_{2},\ldots).
$$

\noindent Hence, we are going to study mainly the strong records in this paper.\\

\noindent First of all, we are going to see that the records times and the record values are Markovian.\\

\bigskip \noindent \textbf{(II) - Markovian properties of the records and the record times}.\\

\noindent We have the following fact concerning the strong record times.

\begin{proposition} \label{lrecMarkov}
The sequence $(U(n))_{n\geq 1}$ of strong record times is a Markovian chain with transition probabilities

\begin{equation*}
p_{t,n}(k,j)=\left\{ 
\begin{tabular}{lll}
$\mathbb{P}\left(X_{j}>X_{k}, \max_{k<h<j} X_{h}\leq X_{k}\right)$ & if & $j>k$, \\ 
 $0$ &  & otherwise. 
\end{tabular}
\right.
\end{equation*}
\end{proposition}

\bigskip \noindent \textbf{Proof}. Let us assume that $n\geq 2$. Let $j\geq n+1$, $1=k_1 < k_2 <\ldots <k_n$. Conditionally on $(U(1)=k_1, \ \ldots ,\ U(n-1)=k_{n-1}, \ U(n)=k_{n})$, the event $(U(n+1)=j)$ depends only on the observations $X_h$, $k_n<h\leq j$ and reduces to

$$
\left(X_{k_{n}+1}\leq X_{k_n},  X_{k_{n}+2}\leq X_{k_n}, \ldots X_{j-1} \leq X_{k_n}, X_{j}>X_{k_n} \right),
$$  

\bigskip \noindent Thus we have

\begin{eqnarray*}
&&\mathbb{P}( U(n+1)=j) / (U(1)=k_1, \ \ldots ,\ U(n-1)=k_{n-1}, \ U(n)=k_{n})\\
&&=\mathbb{P}\left(X_{k_{n}+1}\leq X_{k_n},  X_{k_{n}+2}\leq X_{k_n}, \ldots X_{j-1} \leq X_{k_n}, X_{j}>X_{k_n} \right),
\end{eqnarray*}

\noindent which proves that behavior of $U(n)$ depends only on the most recent past $k_n$, and by the way, provides the probability transition. \\

\noindent Let us move to the record values. We have 

\begin{proposition} \label{lrecMarkov02} 
The sequence $(X^{(n)})_{n\geq 1}$ of strong records is a Markovian chain with transition probability

\begin{equation*}
p_{r,n}(x,A)=\sum_{k=n}^{+\infty} \mathbb{P}\left( (X_{min(j>k, \ X_j>x_n)}  \in A) / (X^{(n)}=x)\right).
\end{equation*}

\noindent where $x$ is a real number and $A$ a Borel set of $\mathbb{R}$.
\end{proposition}

\bigskip \noindent \textbf{Proof}. Let $n\geq 2$. Let $x_1 < x_2 <\ldots <x_n<x$. Conditionally on the intersection $B=(X^{(1)}=x_1, \ \ldots ,\ X^{(n-1)}=x_{n-1}, \ X^{(n)}=x_{n})$, the record $X^{(n+1)}$ is defined on $(U(n)=k)$ with $k\geq n$ by

$$
X^{(n+1)}=X_{min(j>k, \ X_j>x_n)}
$$

\bigskip \noindent so that

$$
\mathbb{P}( (X^{(n+1)}\in A)\cap (U(n)=k)  / B)=\mathbb{P}\left((X_{min(j>k, \ X_j>x_n)})  \in A) \cap (U(n)=k) /(X^{(n)}=x_n) \right).
$$

\bigskip \noindent  We get

$$
\mathbb{P}((X^{(n+1)}\in A)  / B)=\sum_{k=n}^{+\infty} \mathbb{P}( (X_{min(j>k, \ X_j>x_n)})  \in A) \cap (U(n)=k) /(X^{(n)}=x_n)),
$$

\Bin which concludes the proof.\\

\bigskip \noindent Now let us face the general probability laws of the sequences of record values $(X^{(n)})_{n\geq 1}$, of record times
$(U(n))_{n\geq 1}$, of inter-record times $(\Delta(n))_{n\geq 1}=(\Delta_n)_{n\geq 1}$, of the number of record values in a sample $(N(n)+_{n\geq 1}$.\\

\noindent In each case, we give a general probability law regardless to the dependence between the $X_j$'s. Next, we adapt the results to the situation where the $X_j$'s are independent.\\

\noindent Finally, we give detailed results for the \textit{iid} case. In the last case, we apply the previous results but compare our outcomes with formulas with those available in the literature, in particular in books of Ahsanullah.\\

\subsection[Joint Probability Law]{General Joint Cumulative Distribution Functions}

We are going to provide the most general expression, so that we will be able to refine it if we know more about the structure of the dependence of the finite-distributions.

\begin{theorem} \label{FD} For each $n\geq 1$, we have :\\

\noindent (a) The joint \textit{cdf} of the vector of records $X=(X^{(1)}, X^{(2)},\cdots,X^{(n)})^T$ is given, for any $y=(y_1,...,y_n) \in \mathbb{R}^n$,

\begin{eqnarray*}
&&\mathbb{P}(X^{(1)}\leq y_1,\cdots,X^{(n)}\leq y_n) \ \ \ \ \ \ \ \ \ \ \ \ \ \ \ \ \ \ \ \ \ \ \ \ (FD1)\\
&&=\int \mathbb{P}\left(\bigcap_{j=1}^{n}\left(\max_{z_{j-1}+1\leq h \leq z_j} X_h\leq y_j^{\ast}\right)\right) \ d\mathbb{P}_{(U(1),\cdots,U(n))}(z_1,\cdots,z_n),
\end{eqnarray*}

\bigskip \noindent  where $y_i^{\ast}=\wedge_{j=i}^{n} y_j=min(y_i,\cdots,y_n)$, $i \in \{1,\cdots,n\}$.\\

\noindent (b) For $k$-tuple $(n_1,\cdots,n_k=n)$, $1\leq k \leq n$ with $n_0=0<1\leq n_1<\cdots <n_k$, for any $y=(y_1,...,y_k) \in \mathbb{R}^k$, we have

\begin{eqnarray*}
&&\mathbb{P}(X^{(n_1)}\leq y_1,\cdots,X^{(n_k)}\leq y_k) \ \ \ \ \ \ \ \ \ \ \ \ \ \ \ \ \ \ \ \ (FD2)\\
&&=\int \mathbb{P}\left(\bigcap_{j=1}^{k}\left(\max_{z_{j-1}+1\leq h\leq z_j} X_h\leq y_j^{\ast}\right)\right) \ d\mathbb{P}_{(U(n_1),\cdots,U(n_k))}(z_1,\cdots, z_k. \end{eqnarray*}
\end{theorem}

\Bin \textbf{Proof}. Let us prove Formula (b). Let us consider an increasing sequence $(x_j)_{j\geq 1}$. For $1\leq h\leq\ell \leq n$, we define by $A(h,\ell,y)$ the assertion:

\begin{center}
$\biggr($All observations $x_j$, $h\leq j \leq \ell$ is or are less or equal to $y$ $\biggr)$.
\end{center}

\Ni Let us see simple cases. It is not hard to see that for $1\leq n_1 <n_2<\cdots$. We easily see that

$$
(x_1\leq x, x_{n_1}\leq y)=(A(1,1,x) \ and \ A(1, n_1,y))=(A(1,1,min(x,y)) \ and \ A(2, n_1,y))
$$

\Bin and 

\begin{eqnarray*}
&&(x_{n_1}\leq x, x_{n_2}\leq y)=(A(1,n_1,x) \ and \ A(1, n_2,y))=(A(1,n_1,min(x,y))\\
&&and \ A(n_1+1, n_2,y))\\
\end{eqnarray*}

\Bin If we understand the two previous examples, we can see that

\begin{eqnarray} \label{principleRec} 
&&(x_{n_1}\leq x_1, x_{n_2}\leq x_2, \cdots, x_{n_k}\leq x_k)\\
&&=\biggr((A(1,n_1,min(x_1,\cdots,x_k) \ and \ A(n_1+1, n_2,min(x_2,\cdots,x_k) \notag\\
&& and \ A(n_2+1, n_3,min(x_3,\cdots,x_k), \notag\\
&& \cdots, \notag\\
&& A(n_{k-1}+1, n_k,x_k)\biggr). \notag
\end{eqnarray}

\Bin By Applying this simple rule, we have that the event $(Y^{(n_1)}\leq x_1, Y^{(n_2)}\leq x_2, \cdots, Y^{(n_k)}\leq x_k)$ given the event

$$
(U(n_1)=z_1,\cdots,U(n_k)=z_k)
$$

\Bin (with $U(0)=0$ and $z_0=0$) is equivalent to that all the observations in the block of observations from $X_{z_{i}+1}$ to  $X_{z_{i}}$ are less or equal to $y_i^{\ast}=\wedge_{j=i}^{n} y_j$ for $1\leq i \leq n$. We get that $\mathbb{P}(X^{(n_1)}\leq y_1,\cdots,X^{(n_k)}\leq y_k)$ is equal to

$$
\int \mathbb{P}\left(\bigcap_{j=1}^{k}\left(\max_{z_{j-1}+1\leq h \leq z_j} X_h\leq y_j^{\ast}\right)\right) \ d\mathbb{P}_{(U(n_1),\cdots,U(n_k))}(z_1,\cdots,z_k),
$$

\Bin which is Formula (FD2) of which (FD1) is a particular form.\\

\noindent We have the following corollary.

\begin{corollary} \label{FDcoro} We have:\\

\noindent (a) If the random variables $X_j$'s are independent, then for $k$-tuple $(n_1,\cdots,n_k=n)$, $1\leq k \leq n$ with $n_0=0<1\leq n_1<\cdots <n_k$, for any $y=(y_1,...,y_k) \in \mathbb{R}^k$, we have

\begin{eqnarray*}
&&\mathbb{P}(X^{(n_1)}\leq y_1,\cdots,X^{(n_k)}\leq y_k) \ \ \ \ \ \ \ \ \ \ \ \ \ \ \ \ \ \ \ \ (\text{\textit{FDI1}})\\
&&=\int \prod_{j=1}^{k} \prod_{z_{j-1}+1\leq h\leq z_j} F_{X_h}(y_j^{\ast}) \ \ d\mathbb{P}_{(U(n_1),\cdots,U(n_k))}(z_1,\cdots, z_k). 
\end{eqnarray*}

\Bin (b) If the random variables $X_j$'s are \textit{iid} with common \textit{cdf} $F$, then for a $k$-tuple $(n_1,\cdots,n_k=n)$, $1\leq k \leq n$ with $n_0=0<1\leq n_1<\cdots <n_k$, for any $y=(y_1,...,y_k) \in \mathbb{R}^k$, we have

\begin{eqnarray*}
&&\mathbb{P}(X^{(n_1)}\leq y_1,\cdots,X^{(n_k)}\leq y_k) \ \ \ \ \ \ \ \ \ \ \ \ \ \ \ \ \ \ \ \ (\text{\textit{FDI2}})\\
&&=\int \prod_{j=1}^{k} F^{z_{j}-{z_{j-1}}}(y_j^{\ast}) \ \ d\mathbb{P}_{(U(n_1),\cdots,U(n_k))}(z_1,\cdots, z_k). \ \Diamond
\end{eqnarray*}
\end{corollary}

\Bin 

\subsection{Finiteness or Infiniteness of the total number of records} \label{finitenessTNR}

\noindent Let us begin by a general law.

\begin{proposition} \label{nrec01} For each $k\geq 1$, set 

$$
X^{\star}_{k}= \sup_{h>k} X_h.
$$

\Bin and denote

$$
D_{-}=\{(x,y) \in \mathbb{R}^2, \  x\leq y \}.
$$

\Bin We have for any $n\geq 2$,

$$
\mathbb{P}(U(n+1)=+\infty, \ U(n)<\infty)=\sum_{k\geq n} \mathbb{P}_{(X^{\star}_{k},X_k)}(D_{-}) \mathbb{P}(U(n)=k).
$$
\end{proposition}

\noindent \textbf{Proof}. Conditioning on $(U(n)=k)$, $(U(n+1)=+\infty)$ means that all the $X_h$, $h>k$, are less than $X_k$. The proof is ended by the remark

$$
\mathbb{P}(\max_{h>k} X_h\leq X_k)=\mathbb{P}(X^{\star}_{k}\leq X_k)=\mathbb{P}((X^{\star}_{k},X_k) \in D_{-})=\mathbb{P}_{(X^{\star}_{k},X_k)}(D_{-}). \blacksquare
$$

\Bin Let us give an application of Proposition \ref{nrec01} in the independent case.

\begin{proposition} \label{nrec02}
Suppose that $X_1$, $X_2$, \ldots are independent random variables with respective cumulative distribution functions (\textit{cdf}) $F_{j}$, $j\geq 1$. Then, whenever $U(n)$ is finite, we have

$$
\mathbb{P}(U(n+1)=+\infty, \ U(n)<\infty)=\sum_{k\geq n} \left( \int_{\mathbb{R}} \biggr( \prod_{j>k} F_{j}(x) \biggr) d\mathbb{P}_{X_k}(x) \right) \mathbb{P}(U(n)=k).
$$
\end{proposition}

\Bin \textbf{Proof}. Here $X^{\star}_{k}$ and $X_k$ are independent and we have

$$
\mathbb{P}(U(n+1)=+\infty, \ U(n)<\infty)=\sum_{k\geq n} \mathbb{P}_{X^{\star}_{k}} \otimes \mathbb{P}_{X_k}(D_{-}) \mathbb{P}(U(n)=k).
$$

\Bin Let us use Fubini's theorem to have

\begin{eqnarray*}
&&\mathbb{P}_{X^{\star}_{k}} \otimes \mathbb{P}_{X_k}(D_{-})\\
&&=\int_{\mathbb{R}} d\mathbb{P}_{X_k}(x) \int_{\mathbb{R}} 1_{D_{-}}(x,y) d\mathbb{P}_{X^{\star}_{k}}\\
&&=\int_{\mathbb{R}} \mathbb{P}(X^{\star}_{k} \leq x) d\mathbb{P}_{X_k}(x)\\
&&=\int_{\mathbb{R}} \biggr( \prod_{j>k} F_{j}(x) \biggr) d\mathbb{P}_{X_k}(x).
\end{eqnarray*}

\Bin We get the announced result by combining the above lines.\\

\bigskip \noindent Now let us see what happens if the sequence is stationary, that is $F_j=F$ for all $j\geq 1$. Define the lower and the upper endpoints ($\ell$\textit{ep} and \textit{uep}) of $F$ by

$$
uep(F)=\inf \{x \in \mathbb{R}, \ F(x)>1 \} \  and \ uep(F)=\sup \{x \in \mathbb{R}, \ F(x)<1\}
$$

\bigskip \noindent  We have

$$
\int_{\mathbb{R}} \biggr( \prod_{j>k} F_{j}(x) \biggr) d\mathbb{P}_{X_k}(x)=\int_{-\infty}^{uep(F)} F(x)^{+\infty} dF(x).
$$

\bigskip \noindent \noindent But $F(x)^{+\infty}=0$ unless $x=uep(F)$. This gives

$$
\int_{\mathbb{R}} \biggr( \prod_{j>k} F_{j}(x) \biggr) d\mathbb{P}_{X_k}(x)=\int_{-\infty}^{uep(F)} 1_{\{uep(F)\}} dF(x)=\mathbb{P}(X=uep(F)).
$$

\bigskip \noindent  We conclude that

$$
\mathbb{P}(U(n+1)=+\infty, \ U(n)<\infty)=\sum_{k\geq n} \mathbb{P}(X=uep(F)) \mathbb{P}(U(n)=k)=\mathbb{P}(X=uep(F))
$$

\noindent which leads to the simple result:

\begin{proposition} \label{nrec03} 

Suppose that $X_1$, $X_2$, \ldots are independent and identically distributed random variables with common \textit{cdf} $F$ and let $uep(F)$ denote the upper endpoint of $F$. Then

$$
\mathbb{P}(U(n+1)=+\infty, \ U(n)<\infty)=\mathbb{P}(X=uep(F)).
$$

\Bin As a consequence, the sequence of record values (and of record times) is finite if and only if $uep(F)$ is finite and is an atom of $F$, that is $\mathbb{P}_X(uep(F))>0$.
\end{proposition}

\bigskip \noindent \textbf{Consequences}. The number of time records \textit{a.s.} is infinite in the following cases.\\

\noindent (1) $uep(F)=+\infty$.\\

\noindent (2) $uep(F)<+\infty$ but $\mathbb{P}(X=uep(F))=0$. Example : $X \sim \mathcal{U}(0,1)$.\\

\bigskip \noindent The number of time records may be finite in the following cases.\\

\noindent (1) $X$ is discrete and takes a finite number of points.\\

\noindent (2) $X$ is discrete, takes an infinite number of values such the strict values set $\mathcal{V}_X$ of $X$ has a maximum value. We mean by strict values set, the set of points taken by $X$ with a non-zero probability.\\

\newpage
\subsection{Probability law of the sequence of increments of the record times} \label{PLsinterRecordsTimes}
 
\noindent We make the convention that $U(0)=0$. Let $n\geq 2$. If $U(n)$ is finite, we define $\Delta_{n}=U(n)-U(n-1)$. We have:

\begin{proposition} \label{gs_inc01}   If $U(n)$ is finite, then the joint probability law of  

$$
\left(\Delta_{1}, \cdots, \Delta_n\right)
$$

\bigskip \noindent is given by

\begin{eqnarray*}
&&\mathbb{P}(\Delta_{2}=k_2, ..., \Delta_{n}=k_n)\\
&&=\int_{(x_1<x_2<\cdots<x_{n})} \mathbb{P}\left(\left(\bigcap_{1\leq j \leq n-1} \max_{1+\overline{k}_{j}\leq h \leq \overline{k}_{j+1}-1} X_h \leq x_j\right)\right) \ d\mathbb{P}_{(X_{\overline{k}_{1}}, \cdots, X_{\overline{k}_{n}})}(x_1,\cdots,x_{n}),
\end{eqnarray*}

\bigskip \noindent with $k_1=1$, $k_j\geq 1$ for $j \in \{2,\cdots,n\}$ and $\overline{k}_j=k_1+\cdots+k_j$ for $1\leq j \leq n$.
\end{proposition}

\noindent \textbf{Proof}. It is clear that $\Delta_{1}=k_1$ is possible only for $k_1=1$. So in the sequel, we fix $k_1=1$. Let us find the probability law of $(\Delta_{2}, \dots, \Delta_{n})$ through its discrete probability density 
$$
(\Delta_{2}=k_2, ..., \Delta_{n}=k_n).
$$

\bigskip \noindent for $k_j\geq 1$, for $2\leq j \leq n$, $k_1=1$. Let us define 

$$
C_j=\biggr(\max_{1+\overline{k}_{j}\leq h \leq \overline{k}_{j+1}-1} X_h\leq X_{\overline{k}_{j}}, X_{\overline{k}_{j+1}}>X_{\overline{k}_{j}}\biggr), \ j=1,\cdots,n-1
$$

\Bin and

$$
D_j^{\ast}(t)=\left(\max_{1+\overline{k}_{j}\leq h \leq \overline{k}_{j+1}-1} X_h\leq t, \ \ X_{\overline{k}_{j+1}}> t\right), \ t\in \mathbb{R}, \ j=1,\cdots,n-1
$$

\bigskip \noindent We have

\begin{equation}
(\Delta_{2}=k_2, \cdots, \Delta_{n}=k_n)= \bigcap_{1\leq j \leq n-} C_j.
\end{equation}

\bigskip \noindent So, by conditioning by 

$$
Z=(X_{\overline{k}_{1}}, \cdots, X_{\overline{k}_{n}})=(x_1,\cdots,x_{n}), 
$$

\Bin we get

\begin{eqnarray*}
&&\mathbb{P}(\Delta_{2}=k_2, ..., \Delta_{n}=k_n)\\
&&=\int_{(x_1<x_2<\cdots<x_{n})} \mathbb{P}\left(\bigcap_{1\leq j \leq n-1} D_j^{\ast}(x_j)\right) \ d\mathbb{P}_{(X_{\overline{k}_{1}}, \cdots, X_{\overline{k}_{n}})}(x_1,\cdots,x_{n}). \blacksquare
\end{eqnarray*}

\Bin But conditioning by $Z=(X_{\overline{k}_{1}}, \cdots, X_{\overline{k}_{n}})=(x_1,\cdots,x_{n})$, with $x_1<x_2<\cdots<x_n$, we have

$$
D_j^{\ast}(x_j)=\left(\max_{1+\overline{k}_{j}\leq h \leq \overline{k}_{j+1}-1} X_h\leq x_j\right)=:D_j(x_j), \ j=1,\cdots,n-1,
$$

\Bin since each inequality $X_{\overline{k}_{j+1}}>x_j$ are realized by $(x_{j+1}>x_j)$ and this is taken into account in the integration domain. So

\begin{eqnarray*}
&&\mathbb{P}(\Delta_{2}=k_2, ..., \Delta_{n}=k_n)\\
&&=\int_{(x_1<x_2<\cdots<x_{n})} \mathbb{P}\left(\bigcap_{1\leq j \leq n-1} D_j(x_j)\right) \ d\mathbb{P}_{(X_{\overline{k}_{1}}, \cdots, X_{\overline{k}_{n}})}(x_1,\cdots,x_{n}). \blacksquare
\end{eqnarray*}

\bigskip \noindent We have the following corollary in the independent and in the \textit{iid} cases.

\begin{proposition} \label{gs_inc01Ind} Let \bigskip \noindent $k_1=1$, $k_j\geq 1$ for $j \in \{2,\cdots,n\}$ and $\overline{k}_j=k_1+\cdots+k_j$ for $1\leq j \leq n$. Suppose that  $U(n)$ is finite.\\

\noindent (a) If the $X_j$'s are independent, we have

\begin{eqnarray*}
&&\mathbb{P}(\Delta_{2}=k_2, ..., \Delta_{n}=k_n)\\
&&=\int_{(x_1<x_2<\cdots<x_{n})} \biggr\{\prod_{1\leq j\leq n-1} \left(\prod_{1+\overline{k}_{j}\leq h \leq \overline{k}_{j+1}-1} F_{X_h}(x_j)\right)\biggr\}  \  
 d\otimes_{1}^{n}\mathbb{P}_{(X_{\overline{k}_{j}}}x_j).
\end{eqnarray*}

\bigskip \noindent (b) (See \cite{ahsanullah_03}, page 32) If the $X_j$'s are \textit{iid} with common \textit{cdf} $F$, we have for $\overline{k}^{\ast}_{j}=k_2+\cdots+k_j$,

\begin{eqnarray*}
\mathbb{P}(\Delta_{2}=k_2, ..., \Delta_{n}=k_n)=\left(\left(\overline{k}^{\ast}_{n}+1\right) \prod_{2\leq j\leq n} \overline{k}^{\ast}_{j} \right)^{-1} \ \prod_{2\leq j\leq n} 1_{k_j\geq 1}.
\end{eqnarray*}
\end{proposition}

\Bin Before we give the proof, let us introduce the following lemma.

\begin{lemma} \label{lemmaInterRecords} Let us define, for an arbitrary \textit{cdf} $F$ of measure of Lebesgue-Stieljes $\mathbb{L}$ and $k_j \in \mathbb{N}\setminus \{0\}$, for $n\geq 1$,

\begin{eqnarray}
&&\gamma(F,n,k_1,k_2,...., k_n)\\
&&=\int_{(lep(F)\leq x_1<x_2<\cdots<x_{n}\leq uep(F))} \prod_{1\leq j\leq n} F(x_j)^{k_{j}-1} \  d\otimes_{1}^{n}\mathbb{L}^{\otimes n}(x_j).\notag
\end{eqnarray}

\Bin Then 

\begin{equation*}
\gamma(F,n,k_1,k_2,...., k_n)=\left(\prod_{1\leq j \leq n} \overline{k}_j\right)^{-1}. \label{interRecords}
\end{equation*}
\end{lemma}

\noindent \textbf{Proof of Lemma \ref{lemmaInterRecords}}. Let us proceed by induction. For $n=1$ and $k\geq 2$, 
($k_1=1$, $k_2=k$), we clearly have

\begin{eqnarray*}
\gamma(F,n,k)&=&\int_{x<y} F(y)^{k-1} \  d\mathbb{L}^{\otimes 2}(x,y)\\
&=& \int \biggr(\int_{y>x} F(x)^{k-1} \  dF(y)\biggr) \ dF(x)\\
&=& \int \biggr(\left[F(y)^k/k\right]_{x}^{uep(F)}\biggr) \ dF(x)\\
&=& \frac{1}{k}\int \biggr(1-F(x)^k\biggr) \ dF(x)\\
&=&\frac{1}{k} \biggr(1 - \frac{1}{k+1})\\
&=&\frac{1}{k+1}=\frac{1}{\overline{k}_1} \times \frac{1}{\overline{k_2}}.
\end{eqnarray*}

\bigskip \noindent So, Formula \eqref{interRecords} holds for $n=1$. Now, suppose that it holds for $n\geq 1$. Let us prove it for 
$n+1$. We have

\begin{eqnarray*}
&&\gamma(F,n+1,k_1,\cdots, k_{n+1})\\
\end{eqnarray*}

\newpage
\begin{eqnarray*}
&&=\int_{(lep(F)\leq x_1<x_2<\cdots<x_{n+1})\leq uep(F))} \prod_{1\leq j\leq n+1} F(x_j)^{k_{j}-1} 
\ d\otimes_{1}^{n+1}\mathbb{L}(x_j) \\
&&= \int_{(lep(F)\leq x_1<x_2<\cdots<x_{n}\leq uep(F))} \prod_{1\leq j\leq n} F(x_j)^{k_{j}-1} d\otimes_{1}^{n}\mathbb{L}(x_j)\\
&& \times 
\biggr(\int_{(x_{n}\leq x_{n+1}\leq uep(F))}   F(x_{n+1})^{k_{n+1}-1} d\mathbb{L}(x_{n+1}) \biggr)\\
&&= \frac{1}{k_{n+1}} \int_{(lep(F)\leq x_1<x_2<\cdots<x_{n}\leq uep(F))} \prod_{1\leq j\leq n-1} F(x_j)^{k_{j}-1} 
\left(1-F(x_n)^{k_{n+1}}\right) \biggr) \ d\otimes_{1}^{n}\mathbb{L}(x_j) \\
&&= \frac{1}{k_{n+1}}\biggr(\gamma(F,n,k_1,\cdots, k_{n-1}, k_{n})-\gamma(F,n+1,k_1,\cdots, k_{n-1}, k_{n}+k_{n+1}+1)\biggr)\\
&&=\frac{1}{k_{n+1}}\biggr(\prod_{j=1}^{n-1} \overline{k}_j\biggr)^{-1}\biggr(\frac{1}{\overline{k}_n}-\frac{1}{\overline{k}_n+k_{n+1}}\biggr)\\
&&=\biggr(\prod_{j=1}^{n+1} \overline{k}_j\biggr)^{-1}.
\end{eqnarray*}

\Bin The proof of Lemma \ref{lemmaInterRecords} is complete. $\blacksquare$\\

\noindent \textbf{Proof of Proposition of \ref{gs_inc01Ind}}. The first is a translation of the result in Proposition \ref{gs_inc01} when the $X_j$'s are independent. If the $X_j$ are \textit{iid} and if $F$ is continuous, we get, for $\overline{k}^{\ast}_{j}=k_2+\cdots+k_j$, $j\geq 2$, 




\newpage 
\begin{eqnarray*}
&&\mathbb{P}(\Delta_{2}=k_2, ..., \Delta_{n}=k_n)\\
&&= \int_{(x_1<x_2<\cdots<x_{n-1})} \biggr(\prod_{1\leq j\leq n-1} F(x_j)^{k_{j+1}-1}\biggr) \  d\mathbb{P}^{\otimes (n-1)}_{X}(x_j) \int_{x_n>x_{n-1}} dF(x_n)\\
&&+ \int_{(x_1<x_2<\cdots<x_{n-1})} \biggr(\prod_{1\leq j\leq n-1} F(x_j)^{k_{j+1}-1}\biggr)\biggr(1-F(x_{n-1})\biggr) \  d\mathbb{P}^{\otimes (n-1)}_{X}(x_j)\\
&&= \gamma(F,n,k_2,\cdots, k_{n}) - \gamma(F,n,k_2,\cdots, k_{n-1},k_{n}+1)\\
&&=\biggr(\prod_{j=1}^{n-1} \overline{k}_j\biggr)^{-1}\biggr(\frac{1}{\overline{k}_n}-\frac{1}{\overline{k}_{n}+1}\biggr)\\
&&=\biggr((k_n+1) \prod_{j=1}^{n} \overline{k}_j\biggr).
\end{eqnarray*}

\bigskip \noindent $\blacksquare$.\\

\noindent \textbf{Applications}. For $n=2$, we have

\begin{equation}
\mathbb{P}(\Delta_{2}=k)=\frac{1}{k(k+1)}, \label{interRecordsN=1}
\end{equation}

\bigskip \noindent for $n=3$, we get

\begin{equation}
\mathbb{P}\left(\Delta_{2}=k, \Delta_{3}=\ell\right) =\frac{1}{k(k+\ell)(k+\ell+1)}. \label{interRecordsN=2}
\end{equation}

\Bi
\subsection{Probability Law of the sequence of the record times} \label{PLtimesRecords}

\noindent Since we know the probability law of the sequence of the inter-record times, we may get that of the record times by the general formula, with
$\ell_1=1<\ell_2<\cdots<\ell_n$,

\begin{equation} 
\mathbb{P}(U(2)=\ell_2, ..., U(n)=\ell_n) =\mathbb{P}(\Delta_{2}=\ell_2-\ell_1, ..., \Delta_{n}=\ell_n-\ell_{n-1}) \label{lawRtimes}
\end{equation}

\bigskip \noindent Now, we may derive the law of $(U(2), \cdots, U(n))$, $n\geq 2$. Since we have, for $\ell_2<\ell_3<\cdots<\ell_n)$,

\begin{eqnarray} 
\mathbb{P}(U(2)=\ell_2, ..., U(n)=\ell_n) &=&\mathbb{P}(\Delta_{2}=\ell_2-1, ..., \Delta_{n}=\ell_n-\ell_{n-1}) \label{lawRtimes1}.
\end{eqnarray}

\bigskip \noindent By applying Point (b) in proposition \ref{gs_inc01Ind}, we get

\begin{theorem}\label{theolawRtimes} (See \cite{ahsanullah_03}, page 33) If the $X_j$'s are real-valued \textit{iid} random variables, the joint probability law of the record times 

$$
(U(2), \cdots, U(n)),
$$

\Bin for $n\geq 2$, is given by

\begin{eqnarray} 
\mathbb{P}(U(2)=\ell_2, ..., U(n)=\ell_n) = \ell_n^{-1} \prod_{2\leq j \leq n} (\ell_j-1)^{-1} 1_{(2\leq \ell_2<\cdots<\ell_n)}.
\end{eqnarray}
\end{theorem}

\Bin 
\section[Records for real-values \textit{iid} sequences]{Probability laws of strong records from \textit{iid} sequences}

\noindent Here, we are going to give the probability laws of the sequence of record values for independent and identically distributed random variables with common probability law $\mathbb{P}_X$. First, we give the joint cumulative distribution. Secondly, we treat the case where $\mathbb{P}_X$ is absolutely continuous with respect to the Lebesgue measure and finally, the case where $\mathbb{P}_X$ is a discrete probability.\\

\noindent Let us suppose that $X$, $X_1$, $X_2$, ... is a sequence of independent and identically distributed real-valued random variables, defined on the same probability space $(\Omega, \mathcal{A}, \mathbb{P}$), with a common \textsl{cdf} $F$ and common \textit{pdf} $f$.

\subsection[Joint Probability Law]{General Joint Cumulative Distribution Functions}

Before we treat records, let us focus on the sequences of the maxima: $M_n=\max_{1\leq j \leq n} X_j$, $n\geq 1$. We have

\begin{theorem} \label{GRDM1} For each $n\geq 1$, we have :\\

\noindent (a) The joint \textit{cdf} of the vector of records $(M_{1}, M_{2},\cdots,M_{n})^T$ is given, for any $y=(y_1,...,y_n) \in \mathbb{R}^n$, by

$$
\mathbb{P}(M_{1} \leq y_1,\cdots, M_{n}\leq y_n)=\prod_{i}^{n} F\left(\bigwedge_{j=i}^{n} y_j\right). \ \ (\text{\textit{PEX1}})
$$

\bigskip \noindent  where $y^{\ast}_i=\wedge_{j=i}^{n} y_j=min(y_i,\cdots,y_n)$, $1\leq i \leq n$.\\

\noindent (b) For a $k$-tuple $(n_1,\cdots,n_k=n)$, $1\leq k \leq n$ with $n_0=0<1\leq n_1<\cdots <n_k$, for any $y=(y_1,...,y_k) \in \mathbb{R}^k$, we have

$$
\mathbb{P}(M_{n_1}\leq y_1,\cdots, M_{n_k}\leq y_k)=\prod_{j=1}^{k} F^{(n_j-n_{j-1})}\left(\bigwedge_{i=j}^{k}y_i\right). \ \ (\text{\textit{PEX2}})
$$
\end{theorem}

\noindent \textbf{Proof}. We repeat the proof of Theorem \ref{FD} to get (a) by applying the principle described
 in \eqref{principleRec} (Page \pageref{principleRec}). The formula in (b) represents a marginal \textit{cdf} in (\textit{PEX1}). $\square$.\\

\noindent To link this with record values, we notice that $(X^{(1)}, \cdots, X^{(n)})=(M_{U(1)}, \cdots, M_{U(n)})$, for $n\geq 1$. We find again the joint \textsl{cdf} for records values as given in Theorem \ref{FD} in the \textit{iid} case, that is: given $(U(1)=k_1,\cdots,U(n)=k_n)$, we have

$$
(X^{(1)}, \cdots, X^{(n)})=(M_{k_1}, \cdots, M_{k_n})
$$

\bigskip
\noindent and applying (b) in Theorem \ref{GRDM1} allows us to get the conditional law. So we have :

\begin{theorem} \label{GRDMR} Let $n\geq 1$. Define

$$
\Gamma_n=\{(\ell_1,\cdots, \ell_n) \in (\mathbb{N}\setminus \{0\})^n, \ \ \ell_1=1<\ell_2< \cdots < \ell_n\}.
$$

\bigskip \noindent (a) The joint \textit{cdf} of the vector of records $(X^{(1)}, X^{(2)},\cdots,X^{(n)})^T$ is given, for any $y=(y_1,...,y_n) \in \mathbb{R}^n$, by

\begin{eqnarray*}
&&\mathbb{P}(X^{(1)}\leq y_1,\cdots, \ X^{(n)}\leq y_n)\\
&&=\sum_{(\ell_1,\cdots, \ell_n) \in \Gamma_n} \prod_{j=1}^{k} F^{(n_j-n_{j-1})}\left(\bigwedge_{i=j}^{k}y_i\right) 
\mathbb{P}(U(1)=\ell_1,\cdots, U(n)=\ell_n). 
\end{eqnarray*}
\end{theorem}

\subsection[Absolutely continuous records]{Probability laws of strong records from absolutely continuous random variable} $ $\\

\noindent We are going to give the finite distributional probability laws of the sequence of strong records, the individual marginal distributions and different marginal laws involving two or more individual margins.\\

\noindent But before we begin, we wish to explain a general method which will be systematically used. All the computations below which are related to the $n$-th record $X^{(n)}$, $n\geq 2$, are based on conditioning on the past record $X^{(n-1)}$. For this, the reader is supposed to know the general formulas below. Let $A$ be a measurable set of $(\Omega,\mathcal{A})$, $S$ and $T$ be two real valued variables defined on the probability space $(\Omega, \mathcal{A}, \mathbb{P})$. If $T$ has an absolute probability density function $f_T$ with respect to the Lebesgue measure and supported by $\mathcal{V}_X$, we have

\begin{equation}
\mathbb{E}(S)=\int_{\mathcal{V}_X} \mathbb{E}(S|(T=t))  f_T(t) dt \label{fcondAC1}
\end{equation}

\Bin and

\begin{equation}
\mathbb{P}(A)=\int_{\mathcal{V}_X} \mathbb{P}(A|(T=t))  f_T(t) dt. \label{fcondAC2}
\end{equation}

\bigskip \noindent If $T$ is discrete with values set $\mathcal{V}_X=\{x_j, \ j\in J\}$, $J\subset \mathbb{N}$, we have

\begin{equation}
\mathbb{E}(S)=\sum_{j\mathcal{V}_X} \mathbb{E}(S|(T=x_j))  \mathbb{P}(T=x_j) \label{fcondD1}
\end{equation}
 
\noindent and

\begin{equation}
\mathbb{P}(A)=\sum_{j\mathcal{V}_X} \mathbb{P}(A|(T=x_j)) \mathbb{P}(T=x_j). \label{fcondD2}
\end{equation}

\bigskip \noindent For advanced readers, we may use the counting measure $\nu$ on the discrete set $\mathcal{V}_X$ and the \textit{pdf} $f_T(x)=\mathbb{P}(T=x), \ x \in \mathcal{V}_X$, with respect to the counting measure to unify both formulas in 

\begin{equation}
\mathbb{E}(S)=\int_{\mathcal{V}_X} \mathbb{E}(S|(T=t))  f_T(t) d\mu(t) \label{fcondG1}
\end{equation}

\Bin and

\begin{equation}
\mathbb{P}(A)=\int_{\mathcal{V}_X} \mathbb{P}(A|(T=t))  f_T(t) d\mu(t). \label{fcondG2}
\end{equation}
 
\bigskip  \noindent where $\mu$ is the Lebesgue measure ($\mu=\lambda$), If  $T$ has a \textit{pdf} $f_T$ with respect to the Lebesgue measure and supported by $\mathcal{V}_X$, then is the counting measure ($\mu=\nu$) on the discrete values set $\mathcal{V}_T$ if $T$ takes at most countable values.\\

\noindent In order to be in a better position for handling the conditioning on the immediate past $X^{(n-1)}$, we introduce the random variable 

\begin{eqnarray*}
&&N(n-1,n)\\
&&= \text{ Number of observations with time strictly between } U(n-1) \ and  \ U(n).
\end{eqnarray*}

\bigskip \noindent In the sequel, $f$, $\mathcal{V}_X$ and $F$ denote the \textit{pdf} of $X$, its support and its \textit{cdf}, respectively,  in the case of stationary sequences. We also denote

$$
R(x)=-\log(1-F(x)) \  \ and \ \ r(x)=\frac{f(x)}{1-F(x)}, \ x \in \mathcal{V}_X, \ F(x)<1.
$$

\bigskip \noindent where $r(.)$ is the hazard function of $X$, with the obvious relation: 

$$
r(x) \ dx=dR(x), \ x \in \mathcal{V}_X, \ F(x)<1.
$$

\bigskip \noindent Let us give at once the joint absolutely probability law from which marginal probabilities will follow.

\begin{theorem} \label{ADR} 
\noindent (a) The joint probability law of $(X^{(1)}, X^{(2)},\cdots,X^{(n)})^T$, $n\geq 1$, is absolutely continuous to the Lebesgue mesure $\lambda^{\otimes n}$ on $\mathbb{R}^n$ with joint \textit{pdf}, for $y=(y_1,\cdots,y_n)$,

$$
f^{(1,2,\cdots, n)}(y)=\biggr(\prod_{i=1}^{n-1} r(y_j)\biggr) f(y_n) \ 1_{((y_1<y_2<\cdots<y_n)\cap \mathcal{V}_X^n)}(y). \ \ \ (\text{ADR1})
$$ 

\Bin (b) Moreover for any $k$-tuple $(n_1,\cdots,n_k)$, $1\leq k \leq n$ with $1<n_1<\cdots <n_k$, the joint \textit{pdf} with respect to $\lambda^{\otimes k}$ with support $\mathcal{V}_X^k$ is given by 

\begin{eqnarray*}
&&f^{(n_1,n_2,\cdots,n_k)}(y_1,\cdots,y_k)\ \ \  \ \ \ \ \ \ \ \ \ \ \ \ \ \ \ \ \ \ \ \ \ \ \ \ \ \ \ \ \ \ \ \ \ \ \ \ (ADR2)\\
&&=\biggr(\prod_{j=1}^{k-1}r(y_j)\biggr) \biggr(\prod_{j=1}^{k} \frac{(R(y_j)-R(y_{j-1}))^{n_{j}-n_{j-1}-1}}{\Gamma(n_{j}-n_{j-1})}\biggr) f(y_k)\\
&&\times 1_{((y_1<\cdots<y_k)\cap \mathcal{V}_X^k)}(y_1,\cdots, y_k).
\end{eqnarray*}

\Bin with $n_0=0$ and $R(y_0)=0$.\\

\Ni \noindent c) Each $X^{(n)}$, $n\geq 1$, is absolutely continuous with \textit{pdf}

$$
f^{(n)}(x)=\frac{R(x)^{n-1}}{\Gamma(n)} f(x) 1_{\mathcal{V}_X}(x).  \ \ \ (ADR3)
$$ 
\end{theorem}

\noindent \textbf{Remark}. The strict inequality in the domain $(y_1<\cdots<y_k)$ is due to the fact that for an absolutely continuous random vector, the event of the equality of some components is a $null$-event. We may also say that ii comes from the fact that the $y_j$'s are strong records. Actually in the case where the \textit{iid} sequence is associated with an absolutely continuous \textit{pdf}, there are strong records only.\\

\noindent \textbf{Proofs}. Let us organize the proofs into points.\\

\noindent \textbf{Direct Proof of (c)}. This part might be proved as a consequence of (b). But we need to learn how to use the Markov property in a simple case. So, we give a direct proof. Let us begin by the two first cases $n=1$ and $n=2$ and and next, we proceed by induction. For $n=1$, we have that $X^{(1)}=X_1$ and the \textit{pdf} $f$ of $X_1$ is given by
Point $(c)$. For $n=2$, we may use the conditioning formula (\ref{fcondG2}) to have

\begin{eqnarray*}
&&\mathbb{P}(X^{(2)} \in [y-dy/2, y+dy/2])\\
&&=\int_{-\infty}^{y} \mathbb{P}(X^{(2)} \in [y-dy/2, y+dy/2] | X_1=x) f(x) dx\\
&&= \int_{-\infty}^{y} \biggr(\sum_{j=0}^{-\infty} \mathbb{P}(X^{(2)} \in [y-dx/2, y+dy/2], | N(1,2)=j)) | (X_1=x))\\
&& f(x) \biggr) \ dx.\\
\end{eqnarray*}

\Bin The following reasoning will be quoted in the sequel as:\\

\noindent \textbf{(ARG)} \label{arg} The event 
$$
((X^{(2)} \in [y-dx/2, y+dy/2], \ N(1,2)=j)) | (X_1=x))
$$ 

\Bin means exactly that $X_1$ is given and we have $j$ independent random variables $X_h$ less than $x$ and the following random variable (in the enumeration) is in $[y-dx/2, y+dy/2]$. The bounds of the integral come from the fact that the second record in $[y-dy/2, y+dy/2]$ should be greater than $x$ and the fact that no point falls in $[y-dy/2, y+dy/2]$ except the second record based on the continuity assumption. So, the probability of this event is still

$$
F(x)^{j} f(y) dy \ (1+o(1)).
$$

\Bin Then we have, for $0<F(x)<1$

\begin{eqnarray*}
&&\mathbb{P}(X^{(2)} \in [y-dy/2, y+dy/2])\\
&&= \ (1+o(1)) ] \int_{-\infty}^{y} \sum_{j=0}^{-\infty} F(x)^{j} f(y)  f(x) dx dy\\
&&= \ (1+o(1)) \ f(y) dy \int_{-\infty}^{y}  \frac{f(x)}{1-F(x)} dx \\
&&= f(y) dy \ (1+o(1)) \int_{-\infty}^{y}  r(x) dx = f(y) dy \ (1+o(1)) \int_{-\infty}^{y}  dR(x)\\ 
&&=  R(y) f(y)dy \ (1+o(1))
\end{eqnarray*}

\Bin and then

\begin{eqnarray*}
&&\frac{\mathbb{P}(X^{(2)} \in [y-dy/2, y+dy/2])}{dy}=  R(y) f(y) \ (1+o(1)).
\end{eqnarray*}

\Bin We  get the \textit{pdf} of $X^{(2)}$  by letting $dy \rightarrow 0$ and find again Formula (\textit{ADR3}) for $n=2$. Now for the general iteration, we proceed by induction by assuming that (\textit{ADR3}) holds for $n\geq 1$. Let us find the \textit{pdf} for $X^{(n+1)}$. We repeat the same method used for $n=2$ to get

\begin{eqnarray*}
&&\mathbb{P}(X^{(n+1)} \in [y-dy/2, y+dy/2])\\
&&=\int_{-\infty}^{y} \mathbb{P}\biggr(X^{(n+1)} \in [y-dy/2, y+dy/2] \ \biggr| \ X^{(n)}=x\biggr) f^{(n)}(x) dx\\
&&= \int_{-\infty}^{y} \biggr(\sum_{j=0}^{-\infty} \mathbb{P}\biggr(\biggr\{X^{(n+1)} \in [y-dy/2, y+dy/2],\  N(n,n+1)=j\biggr\} \ \biggr| \ X^{(n)}=x)\\
&&\times f^{(n)}(x) \biggr) \ dx.\\
\end{eqnarray*}

\newpage
\noindent By using the argument \textbf{(ARG)} (page \pageref{arg}), we see that the probability of the event 

$$
\biggr(\biggr\{X^{(n+1)} \in [y-dy/2, y+dy/2], \ N(n+1,n)=j\biggr\} \biggr| X^{(n)}=x\biggr)
$$ 

\Bin is 

$$
F(x)^{j} f(y) dy \ (1+o(1)).
$$

\Bin Next, we have 

\begin{eqnarray*}
&&\mathbb{P}(X^{(n+1)} \in [y-dy/2, y+dy/2])\\
&&= \int_{-\infty}^{y-dy/2} \sum_{j=0}^{-\infty} F(x)^{j} f(y) dy \ (1+o(1))  f^{(n)}(x) dx\\
&&= f(y) dy \ (1+o(1))\int_{-\infty}^{y-dy/2}  \frac{f^{(n)}(x)}{1-F(x)} dx\\
&&= f(y) dy \ (1+o(1)) \int_{-\infty}^{y-dy/2}  \frac{r(x)R(x)^{n-1}}{\Gamma(n)} dx\\
&&= (1+o(1)) \ f(y) dy \int_{-\infty}^{y-dy/2}  \frac{d(R(x)^{n})}{n\Gamma(n)}\\ 
&&=  \frac{R(y)^{n}}{\Gamma(n+1)} f(y)dy \ (1+o(1)).
\end{eqnarray*}

\Bin We  get the \textit{pdf} of $X^{(n+1)}$ and find again Formula (\textit{ADR3}) for $n+1$. The proof of (c) by induction is finished.$\square$\\

\noindent Next we have :\\

\noindent \textbf{Proof of (a) : The finite distributional probability law}.\\

\noindent Let us begin by proving the result for $n=2$. We have, for $x$ fixed such that $0\leq F(x)<1$ and $y>x$,

\begin{eqnarray*}
&&\mathbb{P}(X^{(1)} \in [x-dx/2, x+dx/2], X^{(2)} \in [y-dy/2, y+dy/2])\\
&&=\mathbb{P}( (X^{(2)} \in [y-dy/2, y+dy/2]) | (X^{(1)} \in [x-dx/2, x+dx/2]))\\
&& \times \mathbb{P}(X^{(1)} \in [x-dx/2, x+dx/2])\\
&&= \mathbb{P}( (X^{(2)} \in [y-dy/2, y+dy/2]) | (X^{(1)} \in [x-dx/2, x+dx/2]) )\\
&& \times f(x) dx \ (1+o(1)),\\
\end{eqnarray*}

\Bin since $f$ is the \textit{pdf}  of $X^{(1)}$. Using the argument \textbf{(ARG)} (page \pageref{arg}) above, the following set

$$
\biggr(X^{(2)} \in [y-dy/2, y+dy/2]  \biggr| \ X^{(1)} \in [x-dx/2, x+dx/2]\biggr),
$$ 

\bigskip \noindent  when decomposed over the events $N(1,2)=j)$, leads to

\begin{eqnarray*}
&&\mathbb{P}(X^{(1)} \in [x-dx/2, x+dx/2], X^{(2)} \in [y-dy/2, y+dy/2])\\
&&=\sum_{j=0}^{+\infty}\mathbb{P}\biggr( X^{(2)} \in [y-dy/2, y+dy/2], \ N(1,2)=j \ \biggr| \ X^{(1)} \in [x-dx/2, x+dx/2] \biggr)\\
&&\times f(x) dx \ (1+o(1))\\
&&= \sum_{j=0}^{+\infty} F(x-dx/2)^j f(y) dy f(x) dx \ (1+o(1))\\
&&=\frac{f(x)f(y)}{1-F(x+dx/2)} dx dy \ (1+o(1)).\\
\end{eqnarray*}

\Bin As $dx$ and $dy$ go to zero, we exploit the continuity of $F$ at $x$, to have

$$
f^{(1,2)}=r(x)f(y) 1_{(x<y) \cap \mathcal{V}_{X}^{2}},
$$

\bigskip \noindent which proves $(b)$ for $n=2$. Suppose it is true for $n\geq 2$. Let us prove it for $n+1$. Let us fix $x_1<\ldots x_n<x_{n+1}$, all of them in $\mathcal{V}_{X}$ and let us write, for short, $[x_i-dx_i/2, x_i+dx_i/2]=x_{i}\pm dx_{i}/2, \ i=1,\ldots,n$. We have

\begin{eqnarray*}
&&\mathbb{P}(X^{(1)} \in x_{1}\pm dx_{1}/2, \ldots,  X^{(n)} \in x_{n}\pm dx_{n}/2, X^{(n+1)} \in x_{n+1}\pm dx_{n+1}/2)\\
&&=\mathbb{P}\biggr(X^{(n+1)} \in x_{n+1}\pm dx_{n+1}/2) \ \biggr| \ X^{(1)} \in x_{1}\pm dx_{1}/2, \ldots,  X^{(n)} \in x_{n}\pm dx_{n}/2\biggr)\\
&&\times \mathbb{P}(X^{(1)} \in x_{i}\pm dx_{i}/2, \ldots,  X^{(n)} \in x_{n}\pm dx_{n}/2)\\
&&= \mathbb{P}\biggr(X^{(n+1)} \in x_{n+1}\pm dx_{n+1}/2 \ \biggr| \ X^{(1)} \in x_{1}\pm dx_{1}/2, \ldots,  X^{(n)} \in x_{n}\pm dx_{n}/2\biggr) \\
&&\times f^{(1,...,n)}(x_1\cdots,x_n) \ dx_{1} \ldots dx_{n} (1+o(1)).
\end{eqnarray*}

\Bin By the Markov property of the records, we have 

\begin{eqnarray*}
&&\mathbb{P}\biggr(X^{(n+1)} \in x_{n+1}\pm dx_{n+1}/2) \ \biggr| \ X^{(1)} \in x_{1}\pm dx_{1}/2, \ldots,  X^{(n)} \in x_{n}\pm dx_{n}/2\biggr) \\
&&=\mathbb{P}\biggr(X^{(n+1)} \in x_{n+1}\pm dx_{n+1}/2 \ \biggr| \ X^{(n)} \in x_{n}\pm dx_{n}/2\biggr).
\end{eqnarray*}

\bigskip \noindent By using the argument \textbf{(ARG)} (page \pageref{arg}), we obtain 

\begin{eqnarray*}
&&\mathbb{P}(X^{(1)} \in x_{1}\pm dx_{1}/2, \ldots,  X^{(n)} \in x_{n}\pm dx_{n}/2, X^{(n+1)} \in x_{n+1}\pm dx_{n+1}/2)\\
&=&\sum_{j=0}^{+\infty} \mathbb{P}\biggr(X^{(n+1)} \in x_{n+1}\pm dx_{n+1}/2, N(n,n+1)=j \ \biggr| \ X^{(n)} \in x_{n}\pm dx_{n}/2\biggr) \\
&\times& f^{(1,...,n)}(x_1,\ldots,x_n) dx_{1} \ldots dx_{n} (1+o(1))\\
&=& \frac{1+o(1)}{1-F(x_{n}-dx_{n}/2)} f(x_{n+1}) f^{(1,...,n)}(x_1,\ldots,x_n) dx_{1} \ldots dx_{n} dx_{n+1}.
\end{eqnarray*}

\bigskip \noindent The induction hypothesis gives

\begin{eqnarray*}
&&\mathbb{P}(X^{(1)} \in x_{1}\pm dx_{1}/2, \ldots,  X^{(n)} \in x_{n}\pm dx_{n}/2, X^{(n+1)} \in x_{n+1}\pm dx_{n+1}/2)\\
&=& r(x_1)\ldots r(x_{n-1}) \frac{f(x_n) \ (1+o(1))}{1-F(x_{n}-dx_{n}/2)} f(x_{n+1}) dx_{1} \cdots dx_{n+1}.
\end{eqnarray*}

\Bin We get $(a)$ for $n+1$ by letting all the $dx_i$ go to zero for $i=1,...,n$.$\blacksquare$\\

\noindent The combination of Points (a) and (c) allows to have for $y \in \mathcal{V}_X$,

\begin{equation}
\int_{x_1<\ldots<x_{n-1}<y} \left(\prod_{j=1}^{n-1} r(x_i)\right) dx_1\ldots dx_{n-1}=\frac{R(y)^{n-1}}{\Gamma(n)}. \label{gs.20}
\end{equation}

\Bin Indeed, the marginal \textit{pdf} $f^{(n)}$ of $X^{(n)}$ is obtained from the joint \textit{pdf} $f^{(1,\ldots,n)}$ by

\begin{eqnarray*}
f^{(n)}(y)&=&\int_{(x_1<\ldots<x_{n-1}<y)} \left(\prod_{j=1}^{(n-1)} r(x_i)\right) f(y) dx_1 \ldots \ dx_n\\
&=& f(y) \int_{(+\infty<x_1<\ldots<x_{n-1}<y)} \left(\prod_{j=1}^{n-1} r(x_i)\right) dx_1 \ldots \ dx_n.
\end{eqnarray*}

\Bin From Point (c), we make an identification and get Formula (\ref{gs.20}). We may replace $-\infty$ by $z\in \mathcal{V}_X, z<y$ and expect to have

\begin{equation}
\int_{(z<x_1<\ldots<x_{n-1}<y)} \left(\prod_{j=1}^{n-1} r(x_i)\right)\ dx_1 \cdots dx_{n-1}=\frac{(R(y)-R(z))^{n-1}}{\Gamma(n)}. \label{gs.21}
\end{equation}

\bigskip \noindent This is proved in the Appendix Part 2. $\square$\\

\noindent \textbf{Proof of (b) : Distribution of the sub-vector of $(X^{(1)}, X^{(2)},\cdots,X^{(n)})^T$, $n\geq 1$}.\\

\noindent From Point (a), we may find the \textit{pdf} of $(X^{(n_1)}, X^{(n_2)}, \ldots X^{(n_k)})$, denoted $f^{(n_1,n_2,...,n)}$ for $1\leq n_1 <n_2 <\ldots <n_k$. Indeed, we have to integrate the joint
\textbf{pdf} $f^{(1,2,...,n_k)}$ with respect to $dx_i$,
 
$$
i \in \{1,\ldots,n\}\setminus \{n_1,\ldots,n_k\},
$$

\Bin that is, for $x_{n_1}<x_{n_2}<\ldots<x_{n_{k-1}}<x_{n_k}$,

\begin{eqnarray*}
&&f^{(n_1,n_2,...,n_k)}(x_{n_1},x_{n_2},\ldots, x_{n_{k-1}} ,x_{n_k})\\
&&=\int_{x_1<\ldots<x_{n_k}} f^{(1,2,...,n)}(x_1,\ldots,x_{n}) dx_1\ldots dx_{n_{1}-1} \ dx_{n_{1}+1}\\
&&\ldots dx_{n_{2}-1} \ldots dx_{n_{k-1}+1}\ldots dx_{n_{k}-1}\\
&&= r(x_{n_1})\int_{x_1<\ldots<x_{n_{1}-1}<x_{n_{1}}} \prod_{j=1}^{n_{1}-1} r(x_j) \ dx_1 \ldots dx_{n_{1}-1}\\
&&r(x_{n_2})\int_{x_{n_{1}}<x_{n_{1}+1}<\ldots<x_{n_{2}-1}<x_{n_{2}}} \prod_{j=n_1+1}^{n_{2}-1} r(x_j) \ dx_{n_{1}+1} \ldots dx_{n_{2}-1}\\
&&r(x_{n_3})\int_{x_{n_{2}}<x_{n_{2}+1}<\ldots<x_{n_{3}-1}<x_{n_{3}}} \prod_{j=n_2+1}^{n_{3}-1} r(x_j) \ dx_{n_{2}+1}\ldots dx_{n_{3}-1}\\
&& \ldots\\
&&1 \times \int_{x_{n_{k-1}}<x_{n_{k}-1}\ldots<x_{n_{k-1}-1}<x_{n_{k}}} \prod_{j=n_{k-1}+1}^{n_{k}-1} r(x_j) \ dx_{n_{k-1}+1}\ldots 
dx_{n_{k}-1} \ \times f(x_{n_{k}}).\\
\end{eqnarray*}

\bigskip \noindent By applying Formulas (\ref{gs.20}) and (\ref{gs.21}), we get, for $x_{n_0}=lep(F)$,

\begin{eqnarray*}
&&f^{(n_1,n_2,...,n_k)}(x_{n_1},x_{n_2},\ldots, x_{n_{k-1}} ,x_{n_k})\\
&&=f(x_{n_k}) \prod_{j=1}^{k-1} r(x_{n_j})  \prod_{j=1}^{k} \frac{(R(x_{n_{j}})-R(x_{n_{j-1}}))^{n_{j}-n_{j-1}-1}}{\Gamma(n_{j}-n_{j-1})}. \ \ \blacksquare
\end{eqnarray*}

\subsection{Examples of distributions of records}

\noindent A few examples could be useful for studying all the strong records for instance because of the \\

\noindent \textbf{(A) Re-scaling property}. For any sequence $(x_{j})_{1\leq j \leq n}$ of $n$ real numbers, $n\geq 2$, and for any increasing function $h : \mathbb{R} \rightarrow \mathbb{R}$, the record values of $(h(x_{j}))_{1\leq j \leq n}$ are the images by $h$ of the record values of $(x_{j})_{1\leq j \leq n}$ according to the same order and the record times are the same for $(x_{j})_{1\leq j \leq n}$ and $(h(x_{j}))_{1\leq j \leq n}$.\\

\noindent Another interesting rule is the following representation.\\

\noindent \textbf{(B) The Renyi representation}. Let $X$ and $Y$ two real-valued random variables defined on $(\Omega, \mathcal{A}, \mathbb{P})$ with \textit{cdf}'s $F$ and $G$ where $G$ is invertible. Then

$$
X=_d F^{-1}(G(Y)),
$$

\Bin where $F^{-1}(u)=\inf\{x \in [lep(F), uep(F)], \ F(x)\geq u\}$, $u \in ]0,1[$, is the generalized inverse of $F$ and

$$
lep(F)=\inf \{x \in \mathbb{R}, \ F(x)>0\}, \  uep(F)=\sup \{x \in \mathbb{R}, \ F(x)<1\}.
$$

\Bin By applying Point (A) above to random variable we have :

\begin{proposition} \label{repRecGen} Let us consider an \textit{iid} sequence of random variables $(X_{j})_{1\leq j \leq n}$ with common \textit{cdf} $F$, defined on the same probability space $(\Omega, \mathcal{A}, \mathbb{P})$. Let us suppose that $F$ is strictly increasing. Then for any \textit{iid} sequence of random variables $(Y_{j})_{1\leq j \leq n}$ with common \textit{cdf} $G$, we have the following equalities in distribution of the two sequences

\begin{equation}
\{X_1, \cdots, X_n\}=_d \{Y_1, \cdots, Y_n\}, \label{repRecGenRD}
\end{equation}

\Bin and (U(H,j) is the $j$-th record time from the \textit{cdf} $H$),

\begin{equation}
\{U(F,1, \cdots, U(F,n)\}=_d \{U(G,1), \cdots, U(G,n)\}, \label{repRecGenRT}
\end{equation}

\Bin and finally by the re-scaling property for the record values,

\begin{equation}
\{X^{(1)}, \cdots, X^{(n)}\}=_d \{F^{-1}(G(Y^{(1)})), \cdots, F^{-1}(G(Y^{(n)}))\}. \label{repRecGenRV}
\end{equation}
\end{proposition}

\Bin An important example concerns the case where $G$ is the standard exponential \textit{cdf}: $F(x)=1-exp(-x)$, $x\geq 0$. We get

\begin{proposition} \label{repRecExp} Let us consider sequence \textit{iid} sequence of random variables $(X_{j})_{1\leq j \leq n}$ with common \textit{cdf} $F$, defined on the same probability space $(\Omega, \mathcal{A}, \mathbb{P})$. Let us suppose that $F$ is strictly increasing and $G$ is invertible. Then for any \textit{iid} sequence of \textit{iid} standard exponential random variables $(E_{j})_{1\leq j \leq n}$ , we have the equality in distribution of the record times of both sequences

\begin{equation}
\{U(F,1), \cdots, U(F,n)\}=_d \{U(E,1), \cdots, U(F,n)\}, \label{repRecGenRTEXP}
\end{equation}

\Bin and the re-scaling property

\begin{equation}
\{ X^{(1)}, \cdots, X^{(n)}\}=_d \{ F^{-1}(1-exp(-E^{(1)})), \cdots, F^{-1}(1-exp(-E^{(n)}) \}. \label{repRecGenRVEXP}
\end{equation}
\end{proposition}

\Bin Through Formula \eqref{repRecGenRVEXP}, a huge number of problems on records are studied and the Records theory becomes a theoretical study on functions $F^{-1}$ and on exponential random variables.\\

\noindent Now, two interesting and useful examples are given, one for absolutely continuous random variables and one for the discrete random variables: the exponential records (that we are giving right now) and the geometric records to be stated on Page \pageref{recGeom} after the discrete records are presented.\\
 
\noindent \textbf{(C) Exponential records}. Let us consider a sequence of independent standard exponential random variables $(X_{j})_{1\leq j \leq n}$, $n\geq 2$, of intensity $\theta>0$, defined on the same probability space $(\Omega, \mathcal{A}, \mathbb{P})$, By Theorem \ref{ADR} (page \pageref{ADR}) and Formula (\textit{ADR}) therein, we have

\begin{equation}
f_{(X^{(1)}, \cdots, X^{(n)})}(x_1,\cdots, x_n)= \theta^n e^{-\theta x_n} 1_{(0\leq x_1 \leq \cdots \leq x_n)}. \label{pdfRecExpo}
\end{equation}
 
\bigskip
\noindent We see that $(X^{(1)}, \cdots, X^{(n)})$ are the arrival times of a Poisson stochastic process of intensity $\lambda$ and hence inter-arrival times (for $X^{(0)}=0$)

$$
X^{(j)}-X^{(j-1)}, \ j\in\{1,\cdots,n\}, 
$$

\bigskip
\noindent are independent and identically distributed following $\theta$-exponential law.\\

\subsection[Discrete records]{Probability laws of strong records from a discrete random variable}

\vskip 0.5cm
\noindent \textbf{(A) General formula in the iid case}.\\

\noindent Here, we suppose that $X$, $X_1$, $X_2$, ... is a sequence of independent and identically distributed real-valued random variables, defined on the same probability space $(\Omega, \mathcal{A}, \mathbb{P}$), with a common \textit{discrete pdf} $f$ given on the strict support $\mathcal{V}_X=\{x_j, \ j\in J\}\subset \mathcal{R}$, $J\subset \mathbb{N}$, by

\begin{equation}
f(x) = \mathbb{P}(X=x) 1_{\mathcal{V}_X}(x), \ x \in \mathbb{R},
\end{equation}

\bigskip \noindent  with the condition

$$
(\forall x \in \mathcal{V}_X, \ f(x)>0), \ \ (\forall x \notin \mathcal{V}_X, \ f(x)=0) \ \ and \ \ \left(\sum_{x \in \mathcal{V}_X} f(x)=1\right).
$$

\Bin Similarly to the absolutely continuous case, we define for all $x \in \mathbb{R}$,

\begin{eqnarray*}
&&F(x)=\sum_{j \in J, \ x_j\leq x} f(x_j)\\
&&1-F(x)=\sum_{j \in J, \ x_j> x} f(x_j)\\
&&=R(x)=\frac{f(x)}{1-F(x)}, \ \ for \ x \in x \in \mathcal{V}_X \ and \ F(x)<1.
\end{eqnarray*}

\Bin Actually, we mainly use the integration methods, here with respect to the counting measure $\mu$ supported by $\mathcal{V}_X$ and defined by

$$
\mu=\sum_{j \in J} \delta_{x_j} , 
$$

\bigskip \noindent  (where $\delta_{x_j}$ is  the Dirac measure concentrated on $x_j$ with mass one), with respect to which the probability law $\mathbb{P}_X$ has the Radon-Nikodym derivative $f$, that is

$$
d\mathbb{P}_X= f \ d\mu.
$$

\Bin Let us introduce a notation that will replace $R(x)=-\log(1-F(x))$ in the discrete case. We denote for $x<y$,

$$
J(x,y)=]x,y[ \cap \mathcal{V}_X \ \ and \ \ \#J(x,y)=d(x,y).
$$

\Bin For two integers $r<s$, we set

$$
Jo(r,s,x,y)=\{(t_1,\cdots,t_{r-s-1})\in J(x,y)^{r-s-1}, \ t_1<\cdots<t_{r-s-1}\}
$$

\Bin We define $R(x,y)=1$ if $d(x,y)=0$ or $r-s-1=0$. Otherwise, we set 

\begin{equation}
R(r,s,x,y)=\sum_{(t_1,\cdots,t_{r-s-1})\in Jo(r,s,x,y)} \prod_{h=1}^{r-s-1} r(t_h). \ \ \label{Rdiscret}
\end{equation}

\noindent In general, similar rules based only on integration with respect to a measure apply. We get the same formulas but the \textit{pfd}'s are with respect to counting measures. But for pedagogical purposes, redoing the proofs using discrete integration has its own merit. By taking the discrete form of the \textit{pdf}'s in Proposition \ref{ADR}, we have

\begin{proposition} \label{DDR} 

\noindent (a) The joint probability law of $(X^{(1)}, X^{(2)},\cdots,X^{(n)})^T$, $n\geq 1$, has \textit{pdf} with respect to $\mu^{\otimes n}$ with support $\mathcal{V}_X^n$ with discrete \textit{pdf}

$$
f^{(n)}(y_1,y_2,\cdots,y_n)=\biggr(\prod_{i=1}^{n-1} r(y_j)\biggr) f(y_n), 1_{((y_1<y_2<\cdots<y_n)\cap \mathcal{V}_X^n)}(y). \ \ \ (\textit{DDR}_2)
$$ 

\Bin (b) Moreover for any $k$-tuple $(n_1,\cdots,n_k)$, $1\leq k \leq n$ with $1<n_1<\cdots <n_k$, the joint discrete \textit{pdf} with respect to $\mu^{\otimes k}$ with support $\mathcal{V}_X^k$ is given by:

\begin{eqnarray*}
f^{(n_1,n_2,...,n_k)}(y_1,\cdots,y_k)=f(y_n) \prod_{j=1}^{k-1} r(y_{n_j}) \prod_{j=1}^{k} R(n_{j-1},n_{j},y_{n_{j-1}},y_{n_{j}}).
\end{eqnarray*}

\Bin (c) For each $n\geq 1$, the \textit{pdf} with respect to the counting measure supported by $\mathcal{V}_X$ is given by 

$$
f^{(n)}(y)= R(0, n,y_0,y) f(y) 1_{\mathcal{V}_X}(y).  \ \ \ (\textit{DDR}_3)
$$ 
\end{proposition}

\noindent \textbf{Remark}. Here, the strict inequality in the domain $(y_1<\cdots<y_k)$ is due to the fact that the $y_j$'s are strong records. Here we have to distinguish between strong and weak records.\\

\noindent \textbf{Proof}. Let us begin by the joint distribution of $(X^{(1)}, X^{(2)},\cdots,X^{(n)})^T$, $n\geq 1$. For $y_1<y_2<\cdots<y_n$, we set \\

\noindent $U(j,j+1)$ as the number of observations strictly between the $j$-th record $X_{U(j)}$ and the $(j+1)$-th record $X_{U(j+1)}$.\\

\noindent The event

$$
A=(X^{(1)}=y_1, \cdots, X^{(n)}=y_n)
$$

\Bin can be decomposed as

$$
A=\sum_{0\leq p_2<+\infty, \cdots, 0\leq p_{n-1}<+\infty} A \cap (U(j,j+1)=p_j, \ 2\leq j\leq n-1).
$$

\Bin Let us consider the cumulative forms of the $p_j$'s: $p_2^\ast=p_2$, $p_j^\ast=p_2+\cdots+p_j$, $2\leq j \leq n-1$. The event

$$
A \cap (U(j,j+1)=p_j, \ 2\leq j\leq n-1)
$$

\Bin exactly means that the first observation is equal to $y_1$, the $p_2$ next observations are less that $y_1$, the $(p_2^\ast+2)$-th observation is equal to $y_2$, the next $p_3$ observations are less than $y_2$, the $(p_3^\ast+3)$-th, $\cdots$, the $(p_{n-2}^\ast+(n-2))$ is equal to $y_{n-2}$, the $p_{n-1}$ next observations are less than $y_{n-1}$ and finally the $(p_{n-1}^\ast+(n-1))$-th (the $n$-th) observation is equal to $y_n$. So we have

\begin{eqnarray*}
&&\mathbb{P}(A \cap (U(j,j+1)=p_j, \ 2\leq j\ n-1))\\
&&=f(y_1)F(y_1)^{p_2} f(y_2)F(y_2)^{p_3} \cdots f(y_{n-1})F(y_{n-1})^{p_{n-1}} f(y_n).
\end{eqnarray*}

\Bin By summing over the $p_j$'s, we get

$$
\mathbb{P}(A)=\biggr(\frac{f(y_1)}{1-F(f(y_1))} \frac{f(y_2)}{1-F(f(y_2))} \cdots \frac{f(y_{n-1})}{1-F(f(y_{n-1}))}\biggr) f(y_n).
$$ 

\bigskip
\noindent This finishes the proof of Point (1) of the theorem. $\square$\\

\noindent \textbf{Proof of Point (2)}. Let $(n_1,\cdots,n_k)$, $1\leq k \leq n$ with $1<n_1<\cdots <n_k$, be a $k$-tuple and let 

$$
1\leq x_{n_1} < y_{n_2}< \cdots <y_{n_{k-1}} < y_{n_{k}})
$$

\Bin such that $(y_{n_1},\cdots, y_{n_{k}}) \in \mathcal{V}_X^k$. \noindent For 

$$
A=(X^{(n_1)}=x_{n_1},\cdots,X^{(n_k)}=x_{n_k}),
$$

\Bin The marginal distribution of  $(X^{(n_1)}, \cdots,X^{(n_k)})$ is given, for $n_0=0$, $y_{0}=lep(F)$, by

\begin{eqnarray*}
\mathbb{P}(A)&=&\sum_{y_j,\ j\notin \{n_1,\cdots, n_k\}} \prod_{j=1}^{k-1} r(y_j) \ f(y_n) \\
&\times& 1_{(y_1<\cdots <y_{n_1-1}<y_{n_1}<y_{n_1+1}<\cdots <y_{n_2-1}<y_{n_2}<y_{n_2+1}<\cdots<y_{n_k-1}<y_{n_k})}\\
&=& f(y_n) \prod_{j=1}^{k-1} r(y_{n_j}) \prod_{j=1}^{k}\sum_{y_{n_{j-1}+1}<\cdots <y_{n_j-1}} r(y_{n_{j-1}+1})\cdots y_{n_{j}-1}\\
&\times &1_{(y_{n_{j-1}}<y_{n_{j-1}+1}<\cdots <y_{n_j-1}<y_{n_j})}\\
&=&f(y_n) \prod_{j=1}^{k-1} r(y_{n_j}) \prod_{j=1}^{k} R(n_{j-1},n_{j},y_{n_{j-1}},y_{n_{j}}). \ \square
\end{eqnarray*}

\Bin \textbf{Proof of Point (3)}. Let $n\geq 1$. The distribution of $X^{(n)}$ is the $n$-th marginal of the joint distribution $(X^{(1)}, X^{(2)},\cdots,X^{(n)})^T$. So for $n_0=0$, $y_{0}=lep(F)$ and for $y \in \mathcal{V}_X$, we have

\begin{eqnarray*}
\mathbb{P}(X^{(n)}=y)&=&\sum_{y_j,\ 1\leq j \leq n-1\notin \{n_1,\cdots, n_k\}} \prod_{j=1}^{k-1} r(y_j) \ f(y) \\
&\times& 1_{(y_1<\cdots <y_{n-1}<y)}\\
&=& R(0, n,y_0,y) f(y). 
\end{eqnarray*}

\bigskip
\noindent \textbf{(B) Application to geometric records}\label{recGeom}. Let us consider a sequence of \textit{iid} geometric random variables $(X_{j})_{1\leq j \leq n}$, $n\geq 2$, with probability $p=1-q \in ]0,1[$, defined on the same probability space $(\Omega, \mathcal{A}, \mathbb{P})$. Let us recall that each $X_j$ follows the probability law of the number of trials which is necessary to have one success in a Bernoulli trial. Hence the common mass discrete \textit{pdf} $f$ is defined by $f(k)=\mathbb{P}(X_1=k)=q^{k-1} p$, $k\geq 1$. Hence for all $k\geq 1$,

$$
1-F(k)=\sum_{h>k} f(k)= p \sum_{h=k+1}^{+\infty} q^{h-1}= q^{k}
$$

\Bin and

$$
r(k)=\frac{f(k)}{1-F(k)}=p/q.
$$

\Bin Now by applying Theorem \ref{DDR} (Page \pageref{DDR}) and Formula (DDR2) therein, we get for $1\leq k_1 < \cdots < k_n$

$$
f_{(X^{(1)}, \cdots, X^{(n)})}(k_1,\cdots, k_n)= (p/q)^{n-1} q^{k_n-1} p =(p/q)^{n} q^{k_n}
$$

\Bin and thus we have for $1\leq k_1 < \cdots <k_n$,

\begin{equation}
f_{(X^{(1)}, \cdots, X^{(n)})}(k_1,\cdots, k_n)= (p/q) q^{k_n} 1_{(1\leq k_1 < \cdots < k_n)}. \label{pdfRecGeom}
\end{equation}
 
\Bin We also see that $(X^{(1)}, \cdots, X^{(n)})$ are the arrival times of a Bernoulli stochastic process of intensity $\lambda$ and hence inter-arrival times(for $X^{(0)}=0$)

$$
X^{(j)}-X^{(j-1)}, \ j\in\{1,\cdots,n\}, 
$$

\Bin are independent and identically distributed following $p$-geometric laws.\\

\newpage
\section{Records in ordered spaces}

\Ni Let us suppose that we have an ordered probability space $(E, \mathcal{B}, \leq)$. The order relation is denoted by 
$\mathbb{R}$ or by $(\leq)$ and by $x \NC y$, we mean that $x$ and $y$ are not comparable. As well, by \textit{ordered probability space}, we mean that the  $\sigma$-algebra $\mathcal{B}$ is compatible with the partial order relation $(\leq)$ in the following sense: the subsets of $(E, \ \mathcal{B})$ or $(E^k, \mathcal{B}^{\otimes k}$, $k\geq 1$, are measurable:

\begin{eqnarray*}
&&{x} \in \mathcal{B}, \ x \in E\\
&& ]\leftarrow, \ x]=\{y \in E, \ x\leq y \} \in \mathcal{B}, \ x \in E\\
&& ]x, \ \rightarrow ]=\{y \in E, \ y=x \} \in \mathcal{B}, \ x \in E\\
&& N_x=\{y \in E, \ y  \ x \} \in \mathcal{B}, \ x \in E\\
&& N=\{(x,y) \in  E^2, \ x \NC y \} \in \mathcal{B}^{\otimes 2}\\
&& etc.\\ 
\end{eqnarray*}

\Bin For example, this holds with $E=\mathbb{R}^d$, $d\geq 1$, endowed with partial order: 

$$
\biggr(\mathbb{R}^d \ni (x_1,\cdots,x_d) \leq (y_1,\cdots,y_d) \in \mathbb{R}^d\biggr) \Leftrightarrow \biggr(\forall i \in \{1,\cdots,d\}, \ x_i\leq y_i\biggr).
$$

\bigskip \noindent As previously, we work with a sequence of random variables $(Z_n)_{n\geq 1}$ defined on some probability space $(\Omega, \mathcal{A}, \mathbb{P})$ with values in $E$.\\

\subsection{Totally ordered spaces}

\noindent The general characterizations of the probability laws for the record values and the record times in Propositions \ref{lrecMarkov}, \ref{lrecMarkov02} and \ref{gs_inc01} (on pages \pageref{lrecMarkov}, \pageref{lrecMarkov02} and \pageref{gs_inc01} respectively) remain valid.\\

\subsection{Partially ordered spaces}
  
\noindent The definition of record times and record values does not change. But the events based on them will not be as simple as \textsl{less} or \texttt{greater}, since the non comparability will count.\\

\noindent New proofs will not be done again. The results are adapted following the following principles. In a totally ordered space, on $(U(n)=k \ and \ U(n+1)=k+\ell)$, $\ell\geq 1$, we have

\begin{equation}
\forall h \in [k+1, k+\ell-1], \ Z_h\leq Z_k \ and \ Z_{k+\ell}>Z_k. \label{interRecordsTO}
\end{equation}

\Bin But for a partial order, \eqref{interRecordsTO} becomes

\begin{equation}
\forall h \in [k+1, k+\ell-1], \ \biggr( (Z_h\leq Z_k) \cup (Z_h \NC Z_k)\biggr) \ and  \ Z_{k+\ell}>Z_k, \label{interRecordsTG}
\end{equation}

\Bin where the intersection over $h \in [k+1, k+\ell-1]$ is the empty set for $\ell=1$. And we denote

\begin{equation}
C_{k,\ell}= \bigcap_{j=k+1}^{k+\ell} \biggr( (Z_h\leq Z_k) \bigcup (Z_h \NC Z_k)\biggr) \bigcap (Z_{k+\ell}>Z_k). \label{interRecordsTG1}
\end{equation}

\Bin In that regard, Proposition \ref{lrecMarkov} becomes:

\begin{proposition} \label{lrecMarkovG}
The sequence $(U(n))_{n\geq 1}$ of strong record times is a Markovian chain with non-homogeneous transition probabilities

\begin{equation*}
p_{t,n}(k,j)=\left\{ 
\begin{tabular}{lll}
$\mathbb{P}\biggr( \left(Z_{j}>X_{k}\right) \bigcap \bigcap_{k+1\leq h \leq j-1} \left(\left(Z_h \leq Z_k\right) \bigcup \left(Z_h \NC X_k\right)\right)\biggr)$ & if & $j>k$, \\ 
 $0$ &  & otherwise. 
\end{tabular}
\right.
\end{equation*}
\end{proposition}

\bigskip
\noindent The proposition \ref{lrecMarkov02} remains valid as

\begin{proposition} \label{lrecMarkov02G} 
The sequence $(Z^{(n)})_{n\geq 1}$ of strong records is a Markovian chain with transition probabilities

\begin{equation*}
p_{r,n}(x,A)=\sum_{k=n}^{+\infty} \mathbb{P}\left( (Z_{min(j>k, \ Z_j>x_n)}  \in A) / (Z^{(n)}=x)\right),
\end{equation*}

\noindent where $x$ is a real number and $A$ a Borel set of $\mathbb{R}$.
\end{proposition}

\bigskip
\noindent The only concern is related to the computations

\begin{eqnarray}
&&\mathbb{P}\left( (Z_{min(j>k, \ Z_j>x_n)}  \in A) / (Z^{(n)}=x)\right) \label{computeRec}\\
&&=\sum_{\ell=1}^{\infty} \mathbb{P}\biggr((Z_{k+\ell} \in A) \bigcap \bigcap_{1\leq <\ell-1} \biggr((Z_{k+h}\leq Z_k) \bigcup (Z_{k+h}\NC Z_k)\biggr)/ (Z^{(n)}=x)\biggr), \notag
\end{eqnarray} 

\bigskip
\noindent with the convention that: for $\ell=1$ or $\ell=2$, for a collection of sets $C_{\ell}$,

$$
\bigcap_{1\leq <\ell-1} C_{\ell}=\emptyset, 
$$

\Bin \textbf{Remark}. Following \cite{resnickPOS}, the notation

$$
(Z_h>Z_k)^c=\left(\left(Z_h \leq Z_k\right) \bigcup \left(Z_h \NC X_k\right)\right), \ h \geq 1
$$

\bigskip \noindent can be used. \\

\noindent Proposition \ref{gs_inc01} is now

\begin{proposition} \label{gs_inc01G}   If $U(n)$ is finite, then the joint probability law of  

$$
\left(\Delta_{1}, \cdots, \Delta_n\right)
$$

\bigskip \noindent is given by

\begin{eqnarray*}
&&\mathbb{P}(\Delta_{2}=k_2, ..., \Delta_{n}=k_n)\\
&&=\int_{(x_1<x_2<\cdots<x_{n})} 
\mathbb{P}\left(\bigcap_{1\leq j \leq n-1} \bigcap_{1+\overline{k}_j\leq h\leq \overline{k}_{j-1}-1} \biggr( (Z_h\leq Z_k) \bigcup (Z_h \NC Z_k)\biggr)\right)\\
&&d\mathbb{P}_{(X_{\overline{k}_{1}}, \cdots, X_{\overline{k}_{n}})}(x_1,\cdots,x_{n}),
\end{eqnarray*}

\bigskip \noindent with $k_1=1$, $k_j\geq 1$ for $j \in \{2,\cdots,n\}$ and $\overline{k}_j=k_1+\cdots+k_j$ for $1\leq j \leq n$.
\end{proposition}

\bigskip \noindent In such a general case, we cannot go further without knowing the topology of $E$. But the three results will lead to more precise characterizations and more fine description once that topology holds. The first step to take will concern $E=\mathbb{R}^d$, $d>1$.\\

\section{Conclusion}

\noindent This presentation offers a full context of general characterization of the main questions about the probabilistic study of the records and their occurrence times of the sequence of random variables with values in an ordered and measurable space $E$. The literature has a great deal of fine tune results on records on $\mathbb{R}$. We have checked that our record values and record times characterizations remain valid for their counter-parts on $\mathbb{R}$, mostly for independent random variables and \textit{iid} ones. The basis of further and more general results is set. The step to take should concern $E=\mathbb{R}^d$, $d>1$, on the steps of the pioneering works of \cite{resnickPOS}.

\newpage \noindent \textbf{Annexe}.\\

\Bin \textbf{1. Direct proof of Formula \ref{gs.20}}.\\

$$
I_n=\int_{\mathbb{R}^n} \prod_{i=1}^{n} r(x_i) 1_{(x_1<\ldots <x_n<y)}dx_1 \ldots dx_n=\frac{(-\log(1-F(y))^{n}}{\Gamma(n+1)},
$$

\Bin where $R(x)=-\log(1-F(x)), \ x\in \mathbb{R}$, and $\Gamma(n)=(n-1)!$, $n\geq 1$.\\

\noindent Let us show it for $n=1,2,3$. The case $n=1$ is immediate since

\begin{equation}
I_1=\int_{\mathbb{R}} r(x) 1_{(x<y)}dx=\int_{-\infty}^{y} dR(x)=\left[R(x)\right]_{-\infty}^{y}=R(y). \label{recanx01}
\end{equation}

\Bin For $n=2$, we have 

\begin{eqnarray*}
I_{2}&=&\int_{-\infty <x_{1}<x_{2}<y} r(x_1) \left(\int_{x_1}^{y} r(x_2) dx_2\right) dx_1 \notag \\
&=& \int_{-\infty <x_{1}<y} r(x_1) \left(\int_{x_1}^{y} dR(x_2)\right) dx_1 \notag \\
&=& \int_{-\infty <x_{1}<y} r(x_1) (R(y)-R(x_1) dx_1  \ \ (L2)\\
&=& R(y) \int_{-\infty <x_{1}<y} r(x_1) dx_1 -\int_{-\infty <x_{1}<y} r(x_1) R(x_1) dx_1 \notag \\
&=& R(y) \int_{-\infty <x_{1}<y} dR(x_1) -\frac{1}{2}\int_{-\infty <x_{1}<y} d(R(x_1)^2) \notag \\
&=& R(y)^2 - R(y)^2/2=(1/2)R(y)^2. \notag 
\end{eqnarray*}

\Bin For $n=3$, we use the results in Formula (\ref{recanx01}) and Line $(L2)$ of the last blocs of formulas above to get

\begin{eqnarray*}
I_{2}&=&\int_{-\infty <x_1 <x_2<y} r(x_1)r(x_2) \left(R(y)-R(x_2) \right) dx_1 dx_2 \notag \\
&=&\int_{-\infty <x_1 <y} r(x_1) \left( \int_{x_1}^{y} R(y)r(x_2)-r(x_2)R(x_2) \right) dx_1  \notag \\
&=&\int_{-\infty <x_1 <y} r(x_1) \int_{x_1}^{y} \left(R(y)^2-R(y)R(x_1) - (R(y)^2-R(x_1)^2) \right) dx_1  \notag \\
\end{eqnarray*}

\newpage
\begin{eqnarray*}
&=&\int_{-\infty <x_1 <y} \int_{x_1}^{y} \left(R(y)^2r(x_1)-R(y)r(x_1)R(x_1) - (R(y)^2 r(x_1)-r(x_1)R(x_1)^2) \right) dx_1  \notag \\
&=&\int_{-\infty <x_1 <y} \int_{x_1}^{y} \left(R(y)^2 dR(x_1)-R(y)d(R(x_1)^2))/2 - (R(y)^2 dR(x_1)-dR(x_1)^3)/3 \right)  \notag \\
&=& R(y)^3 -(R(y)^3))/2 - R(y)^3+R(y)^3)/3=(1/6)R(y)^3.  \notag \\
\end{eqnarray*}

\Bin  From there, the induction is clear.\\

\bigskip \noindent \textbf{2. Direct proof of Formula (\ref{gs.21}) [see page \pageref{gs.21}]}.\\

\noindent We have to prove that

\begin{equation*}
\int_{z<x_1<\ldots<x_{n-1}<y} \left(\prod_{j=1}^{n-1} r(x_j)\right) dx_1\ldots dx_{n-1}=\frac{(R(y)-R(z))^{n-1}}{\Gamma(n)}.
\end{equation*}

\bigskip \noindent Let us begin with $n=2$, which corresponds to the formula 

\begin{equation*}
\int_{z<x<y} r(x) dx=(R(y)-R(z)),
\end{equation*}

\bigskip \noindent which is obvious since $r(x)=dR(x)/dx$. Suppose it is true for $n\geq 2$. Let us prove it for $n+1$. We have

\begin{eqnarray*}
&&\int_{z<x_1<\ldots<x_{n-1}<x_{n}<y} \left(\prod_{j=1}^{n} r(x_j)\right) dx_1\ldots dx_{n-1} dx_{n}\\
&=&\int_{z<x_1<y} r(x_1) dx_1 \int_{x_1<x_2<\ldots<x_{n}<y} \left(\prod_{j=2}^{n} r(x_j)\right) dx_2\ldots dx_{n}.\\
\end{eqnarray*}

\bigskip \noindent The induction hypothesis gives

$$
\int_{x_1<x_2<\ldots<x_{n}<y} \left(\prod_{j=2}^{n} r(x_i)\right) dx_2\ldots dx_{n}=\frac{(R(y)-R(x_1))^{n-1}}{\Gamma(n)},
$$

\Bin and we get

\begin{eqnarray*}
&&\int_{z<x_1<\ldots<x_{n-1}<x_{n}<y} \left(\prod_{j=1}^{n} r(x_i)\right) dx_1\ldots dx_{n-1} dx_{n}\\
&=&\int_{z<x_1<y} r(x_1) \frac{r(x_1)(R(y)-R(x_1))^{n-1}}{\Gamma(n)} dx_1 .\\
&&\int_{z<x_1<y} r(x_1) \frac{-d(R(y)-R(x_1))^{n}}{n\Gamma(n)}=\frac{(R(y)-R(z))^{n}}{n\Gamma(n+1)}.
\end{eqnarray*}

\bigskip
\noindent The proof of Formula (\ref{gs.21}) is complete.$\blacksquare$\\

\noindent \textbf{Direct proof of Formula (\textit{ADR1}) in Theorem \ref{ADR}, page \pageref{ADR}}.\\

\noindent We define $U(j,j+1)$, for $2\leq j \leq n-1$, as the number of observations strictly between the $j$-th record $X_{U(j)}$ and the $(j+1)$-th record $X_{U(j+1)}$. The event

$$
A=(Y^{(1)} \in y_1\pm dy_1/2, \cdots, Y^{(n)} \in y_n\pm dy_n/2)
$$

\bigskip
\noindent can be decomposed as

$$
A=\sum_{0\leq p_2<+\infty, \cdots, 0\leq p_{n-1}<+\infty} A \cap (U(j,j+1)=p_j, \ 2\leq j\leq n-1).
$$

\Bin Let us accumulate the $p_j$'s as $p_2^\ast=p_2$, $p_j^\ast=p_2+\cdots+p_j$, $2\leq j \leq n-1$. Since the data are continuous, implying that the equality of observations has probability zero, we may and do suppose that each interval $y_j\pm dy_j/2$ contains one observation. So the event

$$
A \cap (U(j,j+1)=p_j, \ 2\leq j\leq n-1)
$$

\bigskip
\noindent exactly means that the first observation is equal to $y_1$, the $p_2$ next observations are less than $y_1$, the $(p_2^\ast+2)$-th observation is in $y_1\pm dy_1/2$, the next $p_3$ are less than $y_2$, the $(p_3^\ast+3)$-th, $\cdots$, the $(p_{n-2}^\ast+(n-2))$ is in $y_2\pm dy_2/2$, the $p_{n-1}$ next observations are less than $y_{n-1}$ and finally the $(p_{n-1}^\ast+(n-1))$-th (the $n$-th) is in 
$y_n\pm dy_n/2$. So we have

\begin{eqnarray*}
&&\mathbb{P}(A \cap (U(j,j+1)=p_j, \ 2\leq j\leq n-1))\\
&&= \mathbb{P}(Y^{(1)}\in y_1\pm dy_1/2) F(y_1)^{p_2} \mathbb{P}(Y^{(2)}\in y_2\pm dy_2/2) F(y_2)^{p_3}\\
 &\cdots& \mathbb{P}(Y^{(n-1)}\in y_{n-1}\pm dy_{n-1}/2) F(y_{n-1})^{p_{n-1}} f(y_n)\mathbb{P}(Y^{(n)}\in y_1\pm dy_n/2. 
\end{eqnarray*}

\Bin By summing over the $p_j$'s first, next by dividing by $(dy_1 \cdots dy_n)$ and finally by letting each $dy_j\rightarrow 0$, we get

$$
\mathbb{P}(A)=\biggr(\frac{f(y_1)}{1-F(f(y_1))} \frac{f(y_2)}{1-F(f(y_2))} \cdots \frac{f(y_{n-1})}{1-F(f(y_{n-1}))}\biggr) f(y_n). \ \blacksquare
$$

\label{fin-art}

\end{document}